# A two-dimensional high-order well-balanced scheme for the shallow water equations with topography and Manning friction

Victor Michel-Dansac[a,b,*], Christophe Berthon[c], Stéphane Clain[d,e], Françoise Foucher[c,f]

[a]*Institut de Mathématiques de Toulouse, Université Toulouse 3 Paul Sabatier, 118 route de Narbonne, 31062 Toulouse Cedex 9, France*
[b]*INSA Toulouse, 135 avenue de Rangueil, 31077 Toulouse Cedex 4, France*
[c]*Laboratoire de Mathématiques Jean Leray, CNRS UMR 6629, Université de Nantes, 2 rue de la Houssinière, BP 92208, 44322 Nantes Cedex 3, France*
[d]*CFUM – Centre of Physics, University of Minho, Azurém Campus, 4800-058 Guimarães, Portugal*
[e]*DMAT – Department of Mathematics, University of Minho, Azurém Campus, 4800-058 Guimarães, Portugal*
[f]*École Centrale de Nantes, 1 rue de La Noë, BP 92101 44321 Nantes Cedex 3, France*

**Abstract**

We develop a two-dimensional high-order numerical scheme that exactly preserves and captures the moving steady states of the shallow water equations with topography or Manning friction. The high-order accuracy relies on a suitable polynomial reconstruction, while the well-balancedness property is based on the first-order scheme from [40, 41], extended to two space dimensions. To get both properties, we use a convex combination between the high-order scheme and the first-order well-balanced scheme. By adequately choosing the convex combination parameter following a very simple steady state detector, we ensure that the resulting scheme is both high-order accurate and well-balanced. The method is then supplemented with a MOOD procedure to eliminate the spurious oscillations coming from the high-order polynomial reconstruction and to guarantee the physical admissibility of the solution. Numerical experiments show that the scheme indeed possesses the claimed properties. The simulation of the 2011 Japan tsunami, on real data, further confirms the relevance of this technique.

## 1. Introduction

This work is concerned with the numerical approximation of the shallow water equations, given in two space dimensions by:

$$\begin{cases} \partial_t h + \boldsymbol{\nabla} \cdot \boldsymbol{q} = 0, \\ \partial_t \boldsymbol{q} + \boldsymbol{\nabla} \cdot \left( \dfrac{\boldsymbol{q} \otimes \boldsymbol{q}}{h} + \dfrac{1}{2} g h^2 \boldsymbol{I} \right) = -g h \boldsymbol{\nabla} Z - k \boldsymbol{q} \|\boldsymbol{q}\| h^{-\eta}, \end{cases} \quad (1.1)$$

where $x$ and $y$ are the space variables, $t$ is the time variable, $h \geq 0$ is the water height, $\boldsymbol{q} = {}^t(q_x, q_y)$ is the water discharge, and $g > 0$ is the gravity constant. The homogeneous system is supplemented with two source terms. The topography source term is given by:

$$\boldsymbol{S}^{\boldsymbol{t}} = -g h \boldsymbol{\nabla} Z, \quad (1.2)$$

---

*Corresponding author
Email addresses:* `victor.michel-dansac@math.univ-toulouse.fr` (Victor Michel-Dansac),
`christophe.berthon@univ-nantes.fr` (Christophe Berthon), `clain@math.uminho.pt` (Stéphane Clain),
`francoise.foucher@ec-nantes.fr` (Françoise Foucher)

where $Z(x)$ is the fixed topography function. In addition, the nonlinear Manning friction source term (see for instance [25]) is given by:
$$\boldsymbol{S^f} = -k\boldsymbol{q}\|\boldsymbol{q}\|h^{-\eta}, \tag{1.3}$$
where $k \geq 0$ is the Manning friction coefficient and $\eta = 7/3$.

To shorten the notations, we rewrite (1.1) under the form
$$\partial_t W + \partial_x F(W) + \partial_y G(W) = \begin{pmatrix} 0 \\ \boldsymbol{S^t}(W) + \boldsymbol{S^f}(W) \end{pmatrix}, \tag{1.4}$$
where the conserved variables $W$ and the physical fluxes $F(W)$ and $G(W)$ are given by:
$$W = \begin{pmatrix} h \\ q_x \\ q_y \end{pmatrix}; \quad F(W) = \begin{pmatrix} q_x \\ \dfrac{q_x^2}{h} + \dfrac{1}{2}gh^2 \\ \dfrac{q_x q_y}{h} \end{pmatrix}; \quad G(W) = \begin{pmatrix} q_y \\ \dfrac{q_x q_y}{h} \\ \dfrac{q_y^2}{h} + \dfrac{1}{2}gh^2 \end{pmatrix}. \tag{1.5}$$

The admissible states space of the 2D shallow water equations is the following convex set:
$$\Omega = \left\{ W = {}^t(h, \boldsymbol{q}) \in \mathbb{R}^3, \ h \geq 0 \right\}, \tag{1.6}$$
which accounts for dry areas when $h = 0$. In addition, we define the water velocity $\boldsymbol{u}$ such that $\boldsymbol{q} = h\boldsymbol{u}$. In dry areas, i.e. when $h$ vanishes, we assume that both $\boldsymbol{u}$ and $\boldsymbol{S^f}$ vanish.

In the context of numerical simulations, the preservation of the steady states of the shallow water equations, obtained by taking a vanishing time derivative in (1.1), is of prime importance. Among these stationary solutions is the well-known lake at rest steady state, which is nothing but a steady solution with a vanishing discharge:
$$\begin{cases} \boldsymbol{q} = \boldsymbol{0}, \\ h + Z = \text{cst}. \end{cases} \tag{1.7}$$

This steady state and its numerical preservation have been the object of much work in the last 25 years, we refer for instance to the non-exhaustive lists [36, 35, 2, 15, 10] in one space dimension and [3, 28, 50, 24] in two space dimensions. The general 2D steady states with nonzero discharge, called moving steady states, are constrained with a vanishing discharge divergence, and their study is quite complex. In this manuscript, we only consider moving steady states in one space dimension.

With notations adapted from (1.1), the 1D shallow water equations read
$$\begin{cases} \partial_t h + \partial_x q = 0, \\ \partial_t q + \partial_x \left( \dfrac{q^2}{h} + \dfrac{1}{2}gh^2 \right) = -gh\partial_x Z - kq|q|h^{-\eta}. \end{cases} \tag{1.8}$$

The 1D moving steady states, obtained by canceling the time derivatives in (1.8), are governed by:
$$\begin{cases} q = q_0, \\ \partial_x \left( \dfrac{q_0^2}{h} + \dfrac{1}{2}gh^2 \right) = -gh\partial_x Z - kq|q|h^{-\eta}, \end{cases} \tag{1.9}$$
where $q_0$, which can be nonzero, is some constant and uniform discharge. Considering a vanishing friction, i.e. setting $k = 0$ in (1.9), we recover the moving steady states with topography, governed by:
$$\begin{cases} q = q_0, \\ \partial_x \left( \dfrac{q_0^2}{h} + \dfrac{1}{2}gh^2 \right) = -gh\partial_x Z. \end{cases} \tag{1.10}$$

If we additionally consider smooth data, the moving topography-only steady states turn out to be governed by the following algebraic relation, which is nothing but a statement of Bernoulli's principle:

$$\frac{q_0^2}{2}\partial_x\left(\frac{1}{h^2}\right) + g\partial_x(h+Z) = 0. \tag{1.11}$$

The numerical preservation of these steady states has also been the object of much work in the last two decades, see for instance [31, 18, 6, 40]. On a flat topography, but with a nonzero friction, we get the friction-only steady states, given by

$$\begin{cases} q = q_0, \\ \partial_x\left(\frac{q_0^2}{h} + \frac{1}{2}gh^2\right) = -kq|q|h^{-\eta}. \end{cases} \tag{1.12}$$

Similarly to the topography-only case, the smoothness assumption allows (1.12) to be rewritten under an algebraic form, as follows:

$$-\frac{q_0^2}{\eta-1}\partial_x h^{\eta-1} + \frac{g}{\eta+2}\partial_x h^{\eta+2} + kq_0|q_0| = 0. \tag{1.13}$$

These steady states are highly nonlinear, and exact preservation is a challenging task, see for instance [5, 41] and references within. Note that, for the case of both source terms, the steady states (1.9) cannot be rewritten under an algebraic form.

In [40, 41], the authors develop a robust numerical scheme able to exactly preserve and capture the smooth steady states associated with the topography and the friction source terms. In addition, this scheme was proven to be entropy-satisfying in [12]. We now aim at providing a high-order extension in two space dimensions, while retaining the robustness property, i.e. the preservation of the water height non-negativity, and the essential well-balancedness property. First steps have been undertaken in [11], where a well-balanced second-order MUSCL extension is proposed. Note that other work has been devoted to the development of high-order schemes which preserve the lake at rest (see for instance [16, 22]) or the moving steady states (see for instance [43, 20, 49, 19]). However, these schemes mostly rely on directly reconstructing the algebraic relation (1.11), which entails the added computational cost of having to solve this nonlinear relation for $h$. In addition, the friction source term is left mostly untreated. Therefore, our goal is to propose a well-balanced high-order strategy, for both friction and topography source terms, that does not rely on solving nonlinear equations, and that is applicable to the shallow water equations for two-dimensional geometries.

The paper is organized as follows. First, in Section 2, after recalling the 1D well-balanced scheme following [40, 41], we design the 2D well-balanced scheme. Then, in Section 3, we discuss the high-order polynomial reconstruction and how to apply it to the numerical scheme. Afterwards, we detail a well-balancedness correction in Section 4, designed so that the resulting scheme is both well-balanced and high-order accurate. Section 5 then presents the MOOD techniques we used to ensure the non-negativity preservation and the elimination of the spurious oscillations caused by the high-order reconstruction. Finally, Section 6 is dedicated to the numerical experiments, designed to test the properties of the scheme, namely its well-balancedness, its high-order accuracy and its robustness. The simulation of the 2011 Japan tsunami, on real data, shows good agreement between the numerical results and the physical measurements.

## 2. A well-balanced scheme in two space dimensions

In this section, we build a two-dimensional (2D) scheme by adapting the one-dimensional (1D) well-balanced scheme from [40, 41]. We recall the 1D scheme for the sake of completeness, and we provide a 2D extension.

### 2.1. Reminder of the 1D well-balanced scheme construction

This section is devoted to recalling the construction of a well-balanced scheme for the 1D shallow water equations with the source terms of topography and friction (1.8) is proposed. The notations we use in this section are derived in a straightforward way from (1.4) – (1.5). First, we briefly recall the framework of 1D Godunov-type schemes. To that end, an approximate Riemann solver made of two intermediate states is introduced. Then, we compute these intermediate states, according to the necessary properties they have to satisfy. The intermediate states are given up to two parameters that represent the topography and friction contributions, whose determination is the focus of the third part of this subsection. Finally, we propose a semi-implicit version of the scheme that allows dealing with transitions between dry and wet areas. These four parts of the construction of the scheme have been detailed in [40, 41], and we sketch here their most important components for the sake of completeness.

*2.1.1. Godunov-type schemes*

Let $\Delta x$ be the uniform space step. We discretize the space domain $\mathbb{R}$ in cells $(x_{i-\frac{1}{2}}, x_{i+\frac{1}{2}})$, of volume $\Delta x$ and of center $x_i$. We consider a piecewise constant approximate solution of (1.8) at time $t = t^n$, denoted by $W_i^n$ in each cell. Thus, approximating solutions of (1.8) on such a mesh amounts to solving a Riemann problem at each interface $x_{i+\frac{1}{2}}$, between the states $W_i^n$ and $W_{i+1}^n$. Exactly solving this Riemann problem is usually not computationally efficient, or even possible. Therefore, we introduce an approximate Riemann solver to compute approximate solutions to this Riemann problem, defined as follows (see Figure 1):

$$\widetilde{W}\left(\frac{x}{t}; W_L, W_R\right) = \begin{cases} W_L & \text{if } x/t < \lambda_L, \\ W_L^* & \text{if } \lambda_L < x/t < 0, \\ W_R^* & \text{if } 0 < x/t < \lambda_R, \\ W_R & \text{if } x/t > \lambda_R, \end{cases} \quad (2.1)$$

where $W_L^*$ and $W_R^*$ are the intermediate states, to be defined later, and $\lambda_L$ and $\lambda_R$ are characteristic velocities, which approximate the hyperbolic wave velocities $u \pm \sqrt{gh}$. To ensure that $\lambda_L$ is negative and $\lambda_R$ is positive, we define the characteristic velocities as follows (see for instance [47]):

$$\begin{aligned} \lambda_L &= \min\left(-|u_L| - \sqrt{gh_L},\ -|u_R| - \sqrt{gh_R},\ -\varepsilon_\lambda\right) \leq -\varepsilon_\lambda, \\ \lambda_R &= \max\left(|u_L| + \sqrt{gh_L},\ |u_R| + \sqrt{gh_R},\ \varepsilon_\lambda\right) \geq -\varepsilon_\lambda, \end{aligned} \quad (2.2)$$

where $\varepsilon_\lambda$ is a positive constant that is taken equal to $10^{-10}$ in the numerical simulations.

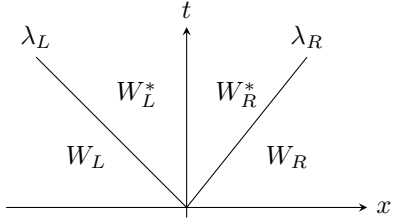

Figure 1: Structure of the approximate Riemann solver.

The piecewise data $W_i^n$ is then evolved according to the approximate Riemann solver, as displayed in Figure 2, leading to a function $W^\Delta$, defined as follows with extended notations:

$$\forall i \in \mathbb{Z}, W^\Delta(x, t^n + t) = \begin{cases} W_{i-\frac{1}{2}}^{R,*} & \text{if } x \in (x_{i-\frac{1}{2}}, x_{i-\frac{1}{2}} + \lambda_{i-\frac{1}{2}}^R t), \\ W_i^n & \text{if } x \in (x_{i-\frac{1}{2}} + \lambda_{i-\frac{1}{2}}^R t, x_{i+\frac{1}{2}} + \lambda_{i+\frac{1}{2}}^L t), \\ W_{i+\frac{1}{2}}^{L,*} & \text{if } x \in (x_{i+\frac{1}{2}} + \lambda_{i+\frac{1}{2}}^L t, x_{i+\frac{1}{2}}). \end{cases} \quad (2.3)$$

Finally, let the time step $\Delta t$ satisfy the following CFL condition:

$$\Delta t \leq \frac{\Delta x}{2\Lambda}, \text{ where } \Lambda = \max_{i \in \mathbb{Z}}(|\lambda_{i+\frac{1}{2}}^L|, \lambda_{i+\frac{1}{2}}^R). \quad (2.4)$$

We define the updated state $W_i^{n+1}$ as the projection of $W^\Delta(x, t^{n+1})$ over the space of piecewise constant functions, namely:

$$W_i^{n+1} := \frac{1}{\Delta x} \int_{x_{i-\frac{1}{2}}}^{x_{i+\frac{1}{2}}} W^\Delta(x, t^{n+1})\, dx.$$

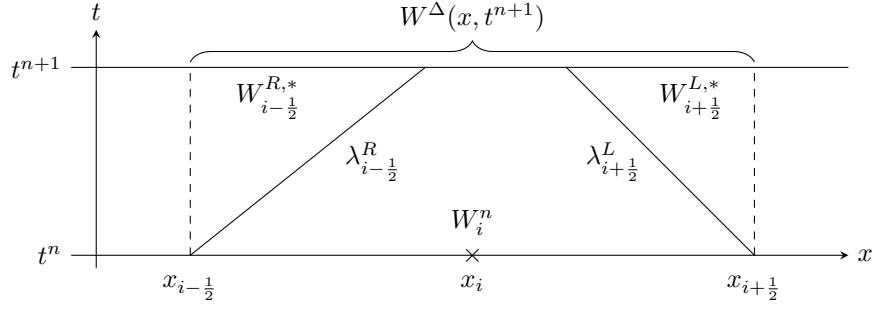

Figure 2: The full Godunov-type scheme using an approximate Riemann solver.

From the definition (2.3) of $W^\Delta$, straightforward computations provide the expression:

$$W_i^{n+1} = W_i^n - \frac{\Delta t}{\Delta x}\left(\lambda_{i+\frac{1}{2}}^L \left(W_{i+\frac{1}{2}}^{L,*} - W_i^n\right) - \lambda_{i-\frac{1}{2}}^R \left(W_{i-\frac{1}{2}}^{R,*} - W_i^n\right)\right). \quad (2.5)$$

To complete the determination of the scheme, one has to give the values of the intermediate states $W_{i+\frac{1}{2}}^{L,*}$ and $W_{i-\frac{1}{2}}^{R,*}$.

### 2.1.2. The intermediate states

According to (2.1), to construct an approximate Riemann solver $\widetilde{W}$ for the Riemann problem between states $W_L$ and $W_R$, we have to determine the two intermediate states, $W_L^*$ and $W_R^*$. We define the intermediate heights and discharges as $W_L^* = {}^t(h_L^*, q_L^*)$ and $W_R^* = {}^t(h_R^*, q_R^*)$ by imposing the three properties the scheme has to satisfy:
- consistency with (1.8);
- well-balancedness;
- non-negativity preservation.

We momentarily assume that $h_L > 0$ and $h_R > 0$.

The consistency of the scheme is obtained by using the consistency relations from [37]. They state that the mean value over a cell of the approximate Riemann solver $\widetilde{W}$, given by (2.1), has to be equal to the average of the exact solution $W_\mathcal{R}$ of the Riemann problem. Therefore, we impose the following relation:

$$\frac{1}{\Delta x}\int_{-\Delta x/2}^{\Delta x/2} \widetilde{W}\left(\frac{x}{\Delta t}; W_L, W_R\right) dx = \frac{1}{\Delta x}\int_{-\Delta x/2}^{\Delta x/2} W_\mathcal{R}\left(\frac{x}{\Delta t}; W_L, W_R\right) dx.$$

Straightforward computations provide the following identities:

$$\begin{aligned}
\lambda_R h_R^* - \lambda_L h_L^* &= (\lambda_R - \lambda_L) h_{HLL}, \\
\lambda_R q_R^* - \lambda_L q_L^* &= (\lambda_R - \lambda_L) q_{HLL} + \frac{1}{\Delta t}\int_0^{\Delta t}\int_{-\Delta x/2}^{\Delta x/2} S(W_\mathcal{R})\, dx\, dt,
\end{aligned} \quad (2.6)$$

where the notation $S$ represents either the topography source term (1.2) or the friction source term (1.3), while ${}^t(h_{HLL}, q_{HLL})$ stands for the intermediate state of the HLL Riemann solver [37]. This intermediate state is defined by

$$\begin{aligned}
(\lambda_R - \lambda_L) h_{HLL} &= \lambda_R h_R - \lambda_L h_L - [q], \\
(\lambda_R - \lambda_L) q_{HLL} &= \lambda_R q_R - \lambda_L q_L - \left[\frac{q^2}{h} + \frac{1}{2}gh^2\right],
\end{aligned}$$

with $[X] = X_R - X_L$ the jump of any quantity $X$.

Equipped with the consistency condition, we now turn to the well-balancedness of the scheme. We observe from the update equation (2.5) that a sufficient condition to get $W_i^{n+1} = W_i^n$, i.e. for the solution to remain stationary, is that $W^{L,*}_{i+\frac{1}{2}} = W_i^n$ and $W^{R,*}_{i-\frac{1}{2}} = W_i^n$ at each interface of the cell of center $x_i$. This situation should be enforced if $W_{i-1}^n$, $W_i^n$ and $W_{i+1}^n$ define a steady state. Let us consider one interface between two states $W_L$ and $W_R$. Here, $W_L$ and $W_R$ are said to define a steady state if the following relations hold:

$$\begin{cases} [q] = 0, \\ \left[\dfrac{q^2}{h} + \dfrac{1}{2}gh^2\right] = \overline{S}, \end{cases} \qquad (2.7)$$

denoting by $\overline{S}$ a consistent average of the source term $S$. We shall determine $\overline{S}$ later for each source term of topography and friction. Note that (2.7) is nothing but a discretization of (1.9), and that $\overline{S}$ stands for the following approximation:

$$\overline{S} \simeq \frac{1}{\Delta x} \frac{1}{\Delta t} \int_0^{\Delta t} \int_{-\Delta x/2}^{\Delta x/2} S(W_{\mathcal{R}})\, dx\, dt.$$

The goal is to derive expressions of $W_L^*$ and $W_R^*$ which ensure that $W_L^* = W_L$ and $W_R^* = W_R$ as soon as $W_L$ and $W_R$ define a steady state according to (2.7).

To that end, after [6, 40, 41], we take $q^* := q_L^* = q_R^*$. Thus, plugging $\overline{S}$ into (2.6) yields:

$$\lambda_R h_R^* - \lambda_L h_L^* = (\lambda_R - \lambda_L) h_{HLL},$$
$$q^* = q_{HLL} + \frac{\overline{S}\Delta x}{\lambda_R - \lambda_L}. \qquad (2.8)$$

An additional relation between the unknowns has to be determined in order to close the system (2.8). According to [40, 41], we choose the following relation:

$$\alpha(h_R^* - h_L^*) = \overline{S}\Delta x, \qquad \text{where} \qquad \alpha = \frac{-\bar{q}^2}{h_L h_R} + \frac{g}{2}(h_L + h_R), \qquad (2.9)$$

with $\bar{q}$ the harmonic mean of $q_L$ and $q_R$:

$$\bar{q} = \begin{cases} \dfrac{2|q_L||q_R|}{|q_L| + |q_R|} \operatorname{sgn}(q_L + q_R) & \text{if } q_L \neq 0 \text{ and } q_R \neq 0, \\ 0 & \text{otherwise.} \end{cases} \qquad (2.10)$$

The relations (2.8) – (2.9) give the intermediate states up to the choice of $\overline{S}$:

$$\begin{aligned} q^* &= q_{HLL} + \frac{\overline{S}\Delta x}{\lambda_R - \lambda_L}, \\ h_L^* &= h_{HLL} - \frac{\lambda_R \overline{S}\Delta x}{\alpha(\lambda_R - \lambda_L)}, \\ h_R^* &= h_{HLL} - \frac{\lambda_L \overline{S}\Delta x}{\alpha(\lambda_R - \lambda_L)}. \end{aligned} \qquad (2.11)$$

Concerning the non-negativity preservation, note that $h_{HLL} > 0$ for the pair $(\lambda_L, \lambda_R)$ defined by (2.2). Now, as presented in [40, 41], we apply the technique introduced in [4, 8] to ensure that $h_L^*$ and $h_R^*$ are non-negative. Indeed, the non-negativity of $h_L^*$ and $h_R^*$ is a sufficient condition to ensure the non-negativity of the updated water height, as evidenced by the expression (2.5) of the scheme. This technique consists in modifying $h_L^*$ and $h_R^*$ in (2.11), to get:

$$\begin{aligned} q^* &= q_{HLL} + \frac{\overline{S}\Delta x}{\lambda_R - \lambda_L}, \\ h_L^* &= \min\left(\left(h_{HLL} - \frac{\lambda_R \overline{S}\Delta x}{\alpha(\lambda_R - \lambda_L)}\right)_+, \left(1 - \frac{\lambda_R}{\lambda_L}\right) h_{HLL}\right), \\ h_R^* &= \min\left(\left(h_{HLL} - \frac{\lambda_L \overline{S}\Delta x}{\alpha(\lambda_R - \lambda_L)}\right)_+, \left(1 - \frac{\lambda_L}{\lambda_R}\right) h_{HLL}\right). \end{aligned} \qquad (2.12)$$

### 2.1.3. Determination of the parameter $\bar{S}$

The intermediate states are characterized by the relations (2.12) up to the parameter $\bar{S}$. We first determine this parameter in the case of the topography source term, and then in the case of the friction source term. We finally show how to combine these two source terms.

*The case of the topography source term.*. We consider two states $W_L$ and $W_R$ defining a steady state with a vanishing friction contribution (i.e. $k = 0$). After (2.7), the following relations govern the steady state:

$$\begin{cases} q_L = q_R = q_0, \\ q_0^2 \left[\dfrac{1}{h}\right] + \dfrac{g}{2}\left[h^2\right] = \bar{S}^t \Delta x, \end{cases} \quad (2.13)$$

where we denote by $\bar{S}^t$ the value of $\bar{S}$ in the case of the topography source term. Note that these relations are a discrete form of (1.10). Moreover, the smooth steady state is also governed by the algebraic relation (1.11), whose discrete form reads:

$$\dfrac{q_0^2}{2}\left[\dfrac{1}{h^2}\right] + g\left[h + Z\right] = 0. \quad (2.14)$$

Equation (2.14) provides an expression of $q_0^2$ that is then plugged into (2.13), and we get:

$$\bar{S}^t \Delta x = -2g[Z]\dfrac{h_L h_R}{h_L + h_R} + \dfrac{g}{2}\dfrac{[h]^3}{h_L + h_R}. \quad (2.15)$$

Note that the expression of $\bar{S}^t$ was previously proposed in the literature (see [6, 7, 40]). However, the average $\bar{S}^t$ has to be consistent with the source term $S^t$, and this property is not verified by the expression (2.15) of $\bar{S}^t$. In order to recover the consistency, we follow the procedure from [6, 7, 40], which consists in introducing a cut-off of the inconsistent term $[h]^3$. Instead of using (2.15), we choose the following expression of $\bar{S}^t$:

$$\bar{S}^t \Delta x := \bar{S}^t(h_L, h_R, Z_L, Z_R, \Delta x)\Delta x = -2g[Z]\dfrac{h_L h_R}{h_L + h_R} + \dfrac{g}{2}\dfrac{[h]_c^3}{h_L + h_R}, \quad (2.16)$$

where $[h]_c$ is the aforementioned cut-off, defined as follows:

$$[h]_c = \begin{cases} h_R - h_L & \text{if } |h_R - h_L| \leq C\Delta x, \\ \operatorname{sgn}(h_R - h_L)C\Delta x & \text{otherwise,} \end{cases} \quad (2.17)$$

with $C$ a positive constant that does not depend on $\Delta x$. Finally, to define the intermediate states when $h_L$ or $h_R$ vanishes, we state again the following result from [40].

**Lemma 1.** *Assume that $W_L$ and $W_R$ define a steady state with a vanishing friction contribution. If $h_L$ or $h_R$ vanishes, we necessarily have $q_0 = 0$. Thus, the quantities $\bar{S}^t$ and $\bar{S}^t/\alpha$ satisfy:*

$$\bar{S}^t \Delta x = -g(Z_R - Z_L)\dfrac{h_L + h_R}{2} \quad \text{and} \quad \dfrac{\bar{S}^t \Delta x}{\alpha} = -(Z_R - Z_L).$$

Note that Lemma 1 only covers the situation of lake at rest steady states governed by (1.7), that is to say where $h_L + Z_L = \text{cst} = h_R + Z_R$. However, according to Figure 3, some physical cases of the lake at rest are not given by (1.7), but governed by

$$\begin{cases} q_0 = h_R = 0, \\ h_L + Z_L \leq Z_R, \end{cases} \quad \text{or} \quad \begin{cases} q_0 = h_L = 0, \\ h_R + Z_R \leq Z_L. \end{cases} \quad (2.18)$$

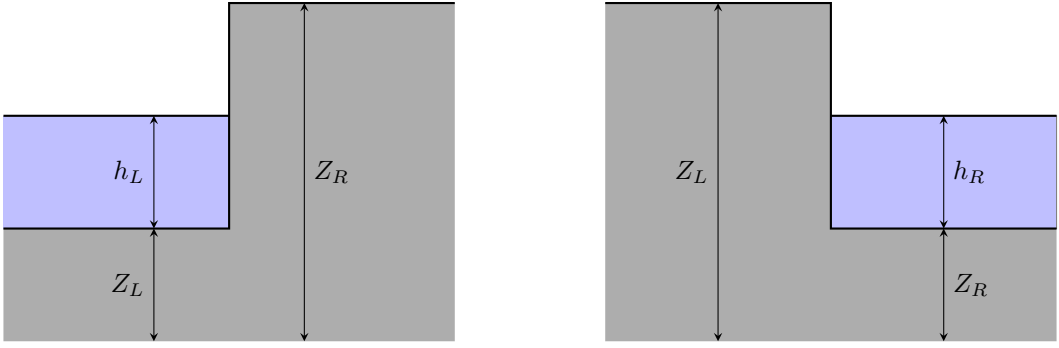

Figure 3: Physical lake at rest configurations not governed by (1.7). Left panel: lake at rest with $h_R = 0$ and $h_L + Z_L \leq Z_R$. Right panel: lake at rest with $h_L = 0$ and $h_R + Z_R \leq Z_L$.

As a consequence, the current expression (2.16) of $\overline{S}^t$ does not provide an exact preservation of the steady states given by (2.18) and depicted in Figure 3. Therefore, we propose an alternative expression of $\overline{S}^t$:

$$\begin{aligned}\overline{S}^t \Delta x &:= \overline{S}^t(h_L, h_R, Z_L, Z_R, \Delta x)\Delta x \\ &= \begin{cases} \dfrac{gh_R^2}{2} & \text{if } q_R = h_L = 0 \text{ and } h_R + Z_R \leq Z_L, \\ -\dfrac{gh_L^2}{2} & \text{if } q_L = h_R = 0 \text{ and } h_L + Z_L \leq Z_R, \\ -g[Z]\dfrac{h_L + h_R}{2} & \text{if } h_L = 0 \text{ or } h_R = 0, \\ -2g[Z]\dfrac{h_L h_R}{h_L + h_R} + \dfrac{g}{2}\dfrac{[h]_c^3}{h_L + h_R} & \text{if } h_L \neq 0 \text{ and } h_R \neq 0. \end{cases}\end{aligned} \quad (2.19)$$

The first two cases of (2.19) are obtained as follows. Assume that $q_L = h_R = 0$ (therefore, $q_R = 0$) and $h_L + Z_L \leq Z_R$. Thus, we take the expression of $\overline{S}^t$ from Lemma 1, and replace $Z_L$ with $\widetilde{Z}_L = Z_R - h_L$, to artificially create a lake at rest configuration governed by (1.7). We immediately obtain the first case of (2.19), and a similar procedure yields the second case. Note that these first two cases indeed exactly preserve the configurations (2.18), since the flux is exactly compensated by $\overline{S}^t$ in (2.13) as soon as any of the two situations (2.18) is considered.

*The case of the friction source term..* We consider two states $W_L$ and $W_R$, with positive water heights, that define a steady state over a flat topography. The approximation $\overline{S}$ is denoted by $\overline{S}^f$ in this friction-only case, and given as follows:

$$\overline{S}^f := \overline{S}^f(h_L, h_R, q_L, q_R) = -k\bar{q}|\bar{q}|\overline{h^{-\eta}}, \quad (2.20)$$

where $\overline{h^{-\eta}}$ is an average depending on $h_L$ and $h_R$, consistent with $h^{-\eta}$ and to be determined, while $\bar{q}$ is the harmonic mean defined by (2.10). Note that $\bar{q} = q_0$ as soon as $W_L$ and $W_R$ define a steady state. Therefore, this steady state is governed by the following relations:

$$\begin{cases} q_L = q_R = q_0, \\ q_0^2\left[\dfrac{1}{h}\right] + \dfrac{g}{2}\left[h^2\right] + kq_0|q_0|\overline{h^{-\eta}}\Delta x = 0, \end{cases} \quad (2.21)$$

which are a discrete version of (1.12). In addition, recall that the algebraic relation (1.13) is satisfied by smooth solutions. The discrete form of this relation reads:

$$-\dfrac{q_0^2}{\eta - 1}\left[h^{\eta - 1}\right] + \dfrac{g}{\eta + 2}\left[h^{\eta + 2}\right] = -kq_0|q_0|\Delta x. \quad (2.22)$$

From (2.22), we obtain an expression of $q_0^2$, which we plug into (2.21), and we obtain:

$$\overline{h^{-\eta}} = \frac{[h^2]}{2} \frac{\eta+2}{[h^{\eta+2}]} - \frac{\bar{\mu}}{k\Delta x} \left( \left[\frac{1}{h}\right] + \frac{[h^2]}{2} \frac{[h^{\eta-1}]}{\eta-1} \frac{\eta+2}{[h^{\eta+2}]} \right),$$

where $\bar{\mu} = \operatorname{sgn}\bar{q}$. To make this expression consistent with $h^{-\eta}$, we once again apply the cutoff (2.17) and we get:

$$\overline{h^{-\eta}} = \frac{[h^2]}{2} \frac{\eta+2}{[h^{\eta+2}]} - \frac{\bar{\mu}}{k\Delta x} [h]_c \left( -\frac{1}{h_L h_R} + \frac{h_L + h_R}{2} \frac{[h^{\eta-1}]}{\eta-1} \frac{\eta+2}{[h^{\eta+2}]} \right), \quad (2.23)$$

The expression (2.23) of $\overline{h^{-\eta}}$ was presented in [41], and proven to be consistent with $h^{-\eta}$. To deal with the case where $h_L$ or $h_R$ vanishes, we follow [41] and make the following statement.

**Statement.** *The quantities $\bar{S}^f$ and $\bar{S}^f/\alpha$ are assumed to vanish as soon as $h_L$ and/or $h_R$ vanish.*

This statement is in agreement with the one made in the introduction that the friction source term vanishes as soon as the water height vanishes.

*The case of both topography and friction source terms..* Thanks to the definitions (2.19) and (2.20) of $\bar{S}^t$ and $\bar{S}^f$, the steady states with topography and friction given by

$$\begin{cases} q_L = q_R = q_0, \\ q_0^2 \left[\dfrac{1}{h}\right] + \dfrac{g}{2}\left[h^2\right] = \bar{S}^t \Delta x + \bar{S}^f \Delta x. \end{cases} \quad (2.24)$$

are exactly preserved by the scheme (2.5) with the intermediate states (2.12). Note that (2.24) is a discretization of the steady relation with topography and friction (1.9). Only this specific discretization will be preserved by the numerical scheme, and as soon as one of the source terms vanish, all steady states for topography or friction, regardless of their discretization, are exactly captured.

Equipped with the two parameters $\bar{S}^t$ and $\bar{S}^f$, we have fully determined the intermediate states (2.12) for $h_L \geq 0$ and $h_R \geq 0$. We thus state the following result, that describes the properties verified by the full scheme.

**Theorem 2.** *Let $W_i^n \in \Omega$ for all $i \in \mathbb{Z}$, with $\Omega$ the admissible states space, given by (1.6). Assume that, for all $i \in \mathbb{Z}$, the intermediate states $W_{i+\frac{1}{2}}^{L,*}$ and $W_{i+\frac{1}{2}}^{R,*}$ satisfy*

$$W_{i+\frac{1}{2}}^{L,*} = \begin{pmatrix} h_L^*\left(W_i^n, W_{i+1}^n\right) \\ q^*\left(W_i^n, W_{i+1}^n\right) \end{pmatrix} \quad \text{and} \quad W_{i+\frac{1}{2}}^{R,*} = \begin{pmatrix} h_R^*\left(W_i^n, W_{i+1}^n\right) \\ q^*\left(W_i^n, W_{i+1}^n\right) \end{pmatrix},$$

*where $h_L^*$, $h_R^*$ and $q^*$ are defined by (2.12). Then the Godunov-type scheme given by (2.5), under the CFL restriction (2.4), satisfies the following properties:*

1. *consistency with the shallow water equations with topography and friction (1.8);*

2. *robustness: $\forall i \in \mathbb{Z}, W_i^{n+1} \in \Omega$;*

3. *well-balancedness: if $(W_i^n)_{i \in \mathbb{Z}}$ defines a steady state according to (2.7), then $\forall i \in \mathbb{Z}, W_i^{n+1} = W_i^n$.*

This result has been proven in [41]. Its proof uses classical ingredients as well as properties gained from the construction of the intermediate states. For the sake of conciseness, we do not recall this proof here.

*2.1.4. Semi-implicitation of the scheme*

The scheme (2.5) – (2.12) is robust, as stated by Theorem 2. However, oscillations due to the stiffness of the source terms occur when simulating transitions between dry and wet areas. To address this issue, we partially follow [41] and introduce an implicit treatment of the friction source term.

*Remark 1.* In [41], both source terms of topography and friction receive an implicit treatment. However, because of the new definition (2.19) of $\bar{S}^t$, such a treatment is not possible anymore for the topography source term. Thus, the three-step scheme from [41] becomes a two-step scheme in the present paper.

We begin by rewriting the scheme (2.5) – (2.12) to exhibit the numerical flux function and the numerical source terms, as follows (see [37] for instance):

$$\begin{pmatrix} h_i^{n+1} \\ q_i^{n+1} \end{pmatrix} = \begin{pmatrix} h_i^n \\ q_i^n \end{pmatrix} - \frac{\Delta t}{\Delta x} \left( \mathcal{F}_{i+\frac{1}{2}}^n - \mathcal{F}_{i-\frac{1}{2}}^n \right) + \Delta t \begin{pmatrix} 0 \\ (S^t)_i^n + (S^f)_i^n \end{pmatrix},$$

where $\mathcal{F}_{i+\frac{1}{2}}^n$ is an approximation of the physical flux at the interface $x_{i+\frac{1}{2}}$, and $(S^t)_i^n$ and $(S^f)_i^n$ are, respectively, approximations of the topography and the friction source terms within the cell $(x_{i-\frac{1}{2}}, x_{i+\frac{1}{2}})$. The numerical flux $\mathcal{F}_{i+\frac{1}{2}}^n$ is defined by:

$$\mathcal{F}_{i+\frac{1}{2}}^n := \mathcal{F}\left(W_i^n, W_{i+1}^n\right) = \frac{1}{2}\left(F(W_i^n) + F(W_{i+1}^n)\right) + \frac{\lambda_{i+\frac{1}{2}}^L}{2}\left(W_{i+\frac{1}{2}}^{L,*} - W_i^n\right) + \frac{\lambda_{i+\frac{1}{2}}^R}{2}\left(W_{i+\frac{1}{2}}^{R,*} - W_{i+1}^n\right). \quad (2.25)$$

In addition, the numerical source terms are defined as follows:

$$(S^t)_i^n = \frac{(\bar{S}^t)_{i-\frac{1}{2}}^n + (\bar{S}^t)_{i+\frac{1}{2}}^n}{2} \quad \text{and} \quad (S^f)_i^n = \frac{(\bar{S}^f)_{i-\frac{1}{2}}^n + (\bar{S}^f)_{i+\frac{1}{2}}^n}{2}, \quad (2.26)$$

where $(\bar{S}^t)_{i+\frac{1}{2}}^n$ and $(\bar{S}^f)_{i+\frac{1}{2}}^n$ are given with clear notations by:

$$(\bar{S}^t)_{i+\frac{1}{2}}^n = \bar{S}^t(h_i^n, h_{i+1}^n, q_i^n, q_{i+1}^n, Z_i, Z_{i+1}, \Delta x),$$
$$(\bar{S}^f)_{i+\frac{1}{2}}^n = \bar{S}^f(h_i^n, h_{i+1}^n, q_i^n, q_{i+1}^n),$$

with the functions $\bar{S}^t$ and $\bar{S}^f$ given by (2.19) and (2.20).

We now introduce a semi-implicit strategy, based on a splitting method (see [15, 47] for more details), which involves an explicit treatment of the hyperbolic part and of the topography source term, and an implicit treatment of the friction source term. Two computational steps are therefore necessary. The first step is devoted to the explicit treatment of the flux and the topography source term, namely:

$$\begin{pmatrix} h_i^{n+\frac{1}{2}} \\ q_i^{n+\frac{1}{2}} \end{pmatrix} = \begin{pmatrix} h_i^n \\ q_i^n \end{pmatrix} - \frac{\Delta t}{\Delta x} \left( \mathcal{F}_{i+\frac{1}{2}}^n - \mathcal{F}_{i-\frac{1}{2}}^n \right) + \Delta t \begin{pmatrix} 0 \\ (S^t)_i^n \end{pmatrix}. \quad (2.27)$$

The second step concerns the friction contribution, which we write as the following initial value problem:

$$\begin{cases} h' = 0, \\ q' = -kq|q|h^{-\eta}, \end{cases} \text{ with initial data } \begin{cases} h(0) = h_i^{n+\frac{1}{2}}, \\ q(0) = q_i^{n+\frac{1}{2}}. \end{cases} \quad (2.28)$$

The initial value problem (2.28) admits an analytic solution, given for all $t \in [0, \Delta t]$ by:

$$\begin{cases} h(t) = h(0), \\ q(t) = \dfrac{h(0)^\eta \, q(0)}{h(0)^\eta + k\, t\, |q(0)|}. \end{cases} \quad (2.29)$$

Evaluating the solution (2.29) at $t = \Delta t$ yields the updated state $W_i^{n+1}$. However, the well-balancedness property on the discharge is lost, since there is no guarantee that $q_i^{n+1} = q_i^n$ when a steady state is assumed. In order to recover this essential behavior, we set:

$$\begin{cases} h_i^{n+1} = h_i^{n+\frac{1}{2}}, \\ q_i^{n+1} = \dfrac{(\bar{h}^\eta)_i^{n+1} \, q_i^{n+\frac{1}{2}}}{(\bar{h}^\eta)_i^{n+1} + k\,\Delta t\, |q_i^{n+\frac{1}{2}}|}, \end{cases} \quad (2.30)$$

where $(\bar{h}^\eta)_i^{n+1}$ is a consistent average of $h_{i-1}^{n+1}$, $h_i^{n+1}$ and $h_{i+1}^{n+1}$ that fulfills the well-balancedness property. More details on the process followed to obtain a relevant average $(\bar{h}^\eta)_i^{n+1}$ are present in [41]. After this step, the quantity $(\bar{h}^\eta)_i^{n+1}$ is defined by:

$$(\bar{h}^\eta)_i^{n+1} = \frac{2k\mu_i^{n+\frac{1}{2}}\Delta x}{k\mu_i^n \Delta x \left(\beta_{i-\frac{1}{2}}^{n+1} + \beta_{i+\frac{1}{2}}^{n+1}\right) - \left(\gamma_{i-\frac{1}{2}}^{n+1} + \gamma_{i+\frac{1}{2}}^{n+1}\right)} + k\Delta t \mu_i^{n+\frac{1}{2}} q_i^n, \quad (2.31)$$

with $\mu_i^n = \mathrm{sgn}(q_i^n)$ and $\mu_i^{n+\frac{1}{2}} = \mathrm{sgn}(q_i^{n+\frac{1}{2}})$, and where $\beta_{i+\frac{1}{2}}^{n+1}$ and $\gamma_{i+\frac{1}{2}}^{n+1}$ are given by:

$$\beta_{i+\frac{1}{2}}^{n+1} = \frac{\eta+2}{2} \frac{\left(h_{i+1}^{n+1}\right)^2 - \left(h_i^{n+1}\right)^2}{\left(h_{i+1}^{n+1}\right)^{\eta+2} - \left(h_i^{n+1}\right)^{\eta+2}},$$

$$\gamma_{i+\frac{1}{2}}^{n+1} = \frac{1}{h_{i+1}^{n+1}} - \frac{1}{h_i^{n+1}} + \beta_{i+\frac{1}{2}}^{n+1} \frac{\left(h_{i+1}^{n+1}\right)^{\eta-1} - \left(h_i^{n+1}\right)^{\eta-1}}{\eta-1}.$$

Equipped with this two-step scheme, we can state the following result.

**Theorem 3.** *The two-step scheme* (2.27) – (2.30) *is consistent with the shallow water equations* (1.8), *robust, and well-balanced.*

A similar result has been proven in [41] for the three-step scheme used in that paper. The proof of this result for the current two-step scheme uses the same ingredients, and we do not present it here for the sake of conciseness.

### 2.2. Two-dimensional extension

We now turn to a two-dimensional extension on a Cartesian grid of the scheme presented in the previous subsection, in order to approximate solutions of (1.1). First, we introduce the notations we use regarding the discretization of the space domain $\mathbb{R}^2$. Then, we present the 2D scheme as a convex combination of 1D schemes.

#### 2.2.1. Space discretization

The discretization of the space domain $\mathbb{R}^2$ consists in a Cartesian mesh of uniform cells, defined by:

$$c_{i,j} = \left(x_{i,j} - \frac{\Delta x}{2}, x_{i,j} + \frac{\Delta x}{2}\right) \times \left(y_{i,j} - \frac{\Delta y}{2}, y_{i,j} + \frac{\Delta y}{2}\right), \quad (2.32)$$

where $(x_{i,j}, y_{i,j})$ is the cell center of $c_{i,j}$. We denote by $|c_{i,j}| = \Delta x \Delta y$ the area of the cell $c_{i,j}$. The piecewise constant approximate solution, within the cell $c_{i,j}$ and at time $t^n$, is denoted by $W_{i,j}^n$.

#### 2.2.2. Construction of the a 2D scheme

We build the two-dimensional extension as a convex combination of one-dimensional schemes. The reader is referred for instance to [44, 10, 14, 13], where such a convex combination is presented for an unstructured mesh. In our particular case of a Cartesian mesh, the first step of the two-step scheme reads as follows:

$$W_{i,j}^{n+\frac{1}{2}} = \frac{1}{4}\left(W_{i+\frac{1}{2},j}^{n+\frac{1}{2}} + W_{i-\frac{1}{2},j}^{n+\frac{1}{2}} + W_{i,j+\frac{1}{2}}^{n+\frac{1}{2}} + W_{i,j-\frac{1}{2}}^{n+\frac{1}{2}}\right), \quad (2.33)$$

with the notations

$$
\begin{aligned}
W_{i-\frac{1}{2},j}^{n+\frac{1}{2}} &= W_{i,j}^n - \frac{4\Delta t}{\Delta x}\left(\mathcal{F}\left(W_{i,j}^n, W_{i,j}^n\right) - \mathcal{F}\left(W_{i,j}^n, W_{i-1,j}^n\right)\right) + 2\Delta t \begin{pmatrix} 0 \\ (\overline{S}_x^t)_{i-\frac{1}{2},j}^n \\ 0 \end{pmatrix}, \\
W_{i+\frac{1}{2},j}^{n+\frac{1}{2}} &= W_{i,j}^n - \frac{4\Delta t}{\Delta x}\left(\mathcal{F}\left(W_{i,j}^n, W_{i+1,j}^n\right) - \mathcal{F}\left(W_{i,j}^n, W_{i,j}^n\right)\right) + 2\Delta t \begin{pmatrix} 0 \\ (\overline{S}_x^t)_{i+\frac{1}{2},j}^n \\ 0 \end{pmatrix}, \\
W_{i,j-\frac{1}{2}}^{n+\frac{1}{2}} &= W_{i,j}^n - \frac{4\Delta t}{\Delta y}\left(\mathcal{G}\left(W_{i,j}^n, W_{i,j}^n\right) - \mathcal{G}\left(W_{i,j}^n, W_{i,j-1}^n\right)\right) + 2\Delta t \begin{pmatrix} 0 \\ 0 \\ (\overline{S}_y^t)_{i,j-\frac{1}{2}}^n \end{pmatrix}, \\
W_{i,j+\frac{1}{2}}^{n+\frac{1}{2}} &= W_{i,j}^n - \frac{4\Delta t}{\Delta y}\left(\mathcal{G}\left(W_{i,j}^n, W_{i,j+1}^n\right) - \mathcal{G}\left(W_{i,j}^n, W_{i,j}^n\right)\right) + 2\Delta t \begin{pmatrix} 0 \\ 0 \\ (\overline{S}_y^t)_{i,j+\frac{1}{2}}^n \end{pmatrix},
\end{aligned}
\quad (2.34)
$$

where $(\overline{S}_x^t)_{i+\frac{1}{2},j}^n$ and $(\overline{S}_y^t)_{i,j+\frac{1}{2}}^n$ are defined by:

$$
\begin{aligned}
(\overline{S}_x^t)_{i+\frac{1}{2},j}^n &= \overline{S}^t\left(h_{i,j}^n, h_{i+1,j}^n, (q_x)_{i,j}^n, (q_x)_{i+1,j}^n, Z_{i,j}, Z_{i+1,j}, \Delta x\right), \\
(\overline{S}_y^t)_{i,j+\frac{1}{2}}^n &= \overline{S}^t\left(h_{i,j}^n, h_{i,j+1}^n, (q_y)_{i,j}^n, (q_y)_{i,j+1}^n, Z_{i,j}, Z_{i,j+1}, \Delta y\right),
\end{aligned}
$$

with $\overline{S}^t$ defined by (2.19). Note that (2.34) represents a collection of four one-dimensional schemes, and that (2.33) is nothing but a convex combination of these schemes.

In (2.34), the numerical flux function $\mathcal{F}$ is the 1D function in the $x$-direction defined in (2.25) and the numerical flux in the $y$-direction, $\mathcal{G}$, derives from $\mathcal{F}$ by using classical rotational invariance properties (see for instance [30]). The equation (2.33) can be rewritten using (2.34) under the following classical form:

$$
\begin{aligned}
W_{i,j}^{n+\frac{1}{2}} = W_{i,j}^n &- \frac{\Delta t}{\Delta x}\left(\mathcal{F}\left(W_{i,j}^n, W_{i+1,j}^n\right) - \mathcal{F}\left(W_{i-1,j}^n, W_{i,j}^n\right)\right) \\
&- \frac{\Delta t}{\Delta y}\left(\mathcal{G}\left(W_{i,j}^n, W_{i,j+1}^n\right) - \mathcal{G}\left(W_{i,j-1}^n, W_{i,j}^n\right)\right) \\
&+ \Delta t \begin{pmatrix} 0 \\ (\boldsymbol{S}^t)_{i,j}^n \end{pmatrix},
\end{aligned}
\quad (2.35)
$$

together with

$$
(\boldsymbol{S}^t)_{i,j}^n = \frac{1}{2}\begin{pmatrix} (\overline{S}_x^t)_{i-\frac{1}{2},j}^n + (\overline{S}_x^t)_{i+\frac{1}{2},j}^n \\ (\overline{S}_y^t)_{i,j-\frac{1}{2}}^n + (\overline{S}_y^t)_{i,j+\frac{1}{2}}^n \end{pmatrix}.
$$

The time step is chosen to satisfy the following CFL condition:

$$
\Delta t = \frac{\delta}{2\Lambda}, \quad (2.36)
$$

where $\delta = \min(\Delta x, \Delta y)$ and $\Lambda$ stands for the maximum of all characteristic velocities at each interface.

To build the second step of the scheme, we address the implicitation of the friction contribution. Similarly to the 1D case, we solve the following initial value problem:

$$
\begin{cases} h' = 0, \\ \boldsymbol{q}' = -k\,\boldsymbol{q}\|\boldsymbol{q}\|h^{-\eta}, \end{cases} \text{with initial data } \begin{cases} h(0) = h_{i,j}^{n+\frac{1}{2}}, \\ \boldsymbol{q}(0) = \boldsymbol{q}_{i,j}^{n+\frac{1}{2}}. \end{cases}
$$

This initial value problem again admits an analytic solution, given for $t \in [0, \Delta t]$ by:

$$
\begin{cases} h(t) = h(0), \\ \boldsymbol{q}(t) = \dfrac{h(0)^\eta\,\boldsymbol{q}(0)}{h(0)^\eta + k\,t\,\|\boldsymbol{q}(0)\|}. \end{cases} \quad (2.37)
$$

We slightly modify the expression of the discharge given by (2.37) to recover the well-balancedness property, and, the updated state reads, with $\boldsymbol{q} = {}^t(q_x, q_y)$:

$$\begin{cases} h_{i,j}^{n+1} = h_{i,j}^{n+\frac{1}{2}}, \\ (q_x)_{i,j}^{n+1} = \dfrac{(\overline{h_x^\eta})_{i,j}^{n+1} \, (q_x)_{i,j}^{n+\frac{1}{2}}}{(\overline{h_x^\eta})_{i,j}^{n+1} + k\,\Delta t \, \left\|\boldsymbol{q}_{i,j}^{n+\frac{1}{2}}\right\|}, \\ (q_y)_{i,j}^{n+1} = \dfrac{(\overline{h_y^\eta})_{i,j}^{n+1} \, (q_y)_{i,j}^{n+\frac{1}{2}}}{(\overline{h_y^\eta})_{i,j}^{n+1} + k\,\Delta t \, \left\|\boldsymbol{q}_{i,j}^{n+\frac{1}{2}}\right\|}, \end{cases} \quad (2.38)$$

with $(\overline{h_x^\eta})_{i,j}^{n+1}$ and $(\overline{h_y^\eta})_{i,j}^{n+1}$ given by (2.31) in the $x$- and $y$-direction, respectively.

The 2D scheme is now complete. To state its properties, we require the following definition.

**Definition.** The vector $(W_{i,j}^n)_{(i,j) \in \mathbb{Z}^2}$ is said to define a steady state in the $x$-direction if:

- $\forall (i,j) \in \mathbb{Z}^2$, $W_{i,j+1}^n = W_{i,j}^n$;
- $\forall (i,j) \in \mathbb{Z}^2$, $(q_y)_{i,j}^n = 0$;
- $\forall (i,j) \in \mathbb{Z}^2$, the pairs ${}^t(h_{i,j}^n, (q_x)_{i,j}^n, Z_{i,j})$ and ${}^t(h_{i+1,j}^n, (q_x)_{i+1,j}^n, Z_{i+1,j})$ satisfy (2.24) or Figure 3.

Similarly, $(W_{i,j}^n)_{(i,j) \in \mathbb{Z}^2}$ is said to define a steady state in the $y$-direction if:

- $\forall (i,j) \in \mathbb{Z}^2$, $W_{i+1,j}^n = W_{i,j}^n$;
- $\forall (i,j) \in \mathbb{Z}^2$, $(q_x)_{i,j}^n = 0$;
- $\forall (i,j) \in \mathbb{Z}^2$, the pairs ${}^t(h_{i,j}^n, (q_y)_{i,j}^n, Z_{i,j})$ and ${}^t(h_{i,j+1}^n, (q_y)_{i,j+1}^n, Z_{i,j+1})$ satisfy (2.24) or Figure 3.

The following result is then satisfied by the 2D scheme:

**Theorem 4.** *Under the CFL condition (2.36), the following properties are satisfied by the two-dimensional two-step scheme (2.35) – (2.38).*

1. *Robustness: if $W_{i,j}^n \in \Omega$ for all $(i,j) \in \mathbb{Z}^2$, then $W_{i,j}^{n+1} \in \Omega$ for all $(i,j) \in \mathbb{Z}^2$.*
2. *Well-balancedness by direction: if $(W_{i,j}^n)_{(i,j) \in \mathbb{Z}^2}$ defines a steady state in the $x$- or the $y$-direction, then for all $(i,j) \in \mathbb{Z}^2$, $W_{i,j}^{n+1} = W_{i,j}^n$.*

*Proof.* The proof of this theorem revolves rewriting the first step (2.35) of the 2D scheme under the form of a convex combination of 1D schemes, given by (2.33) and (2.34). Then, each 1D scheme defined by (2.34) enjoys the same properties as the truly 1D scheme (2.27). In addition, these properties are satisfied by the convex combination (2.33). Therefore, the robustness of the two-step scheme is immediate, since the 1D schemes are robust.

In order to establish the well-balancedness, assume that $(W_{i,j}^n)_{(i,j) \in \mathbb{Z}^2}$ defines a steady state in the $x$-direction. Therefore, the sum of the vertical fluxes in (2.35) and the $y$ contribution of the topography vanish. Thus, the first step of the scheme turns out to be the 1D first step following the $x$-direction. Then, the $y$ contribution of the friction source term vanishes, leaving only the $x$ contribution, which is the same as in the 1D case. Therefore, Theorem 3 applies, and $W_{i,j}^{n+1} = W_{i,j}^n$ for all $(i,j) \in \mathbb{Z}^2$. A similar chain of arguments can be applied to prove the preservation of the steady states in the $y$-direction. $\square$

*Remark 2.* Note that the lake at rest steady state (1.7) is a particular case the steady states defined by direction. Indeed, taking a steady state in the $x$-direction with $q_x = 0$ yields the lake at rest. Similarly, the lake at rest is obtained by taking a steady state in the $y$-direction with $q_y = 0$. Therefore, from Theorem 4, any two-dimensional lake at rest steady state given by (1.7) or by Figure 2.18 is exactly preserved by the 2D scheme.

*Remark 3.* The scheme we have built preserves the 1D moving steady states, in addition to the 2D steady states at rest. The case of the fully 2D steady states, where the discharge is divergence-free instead of being merely constant, is a much more arduous task, and lies outside the scope of this manuscript.

## 3. High-order strategy

Equipped with the first-order 2D scheme built in the previous section, we now turn to building a high-order 2D scheme. First, we focus on the reconstruction strategy used to obtain reconstructed variables. Next, we present the high-order scheme that takes advantage of these reconstructed variables.

From now on and until the end of the paper, the degree of the polynomial reconstruction is denoted by $d \geq 0$. Moreover, for the sake of simplicity, we denote the reconstructed variables by $\varphi \in \{h, q_x, q_y, h+Z\}$. Note that a polynomial reconstruction of degree $d$ provides a scheme of order $(d+1)$.

### 3.1. Obtaining reconstructed variables

We begin by presenting the polynomial reconstruction strategy, introduced in [21, 26] (see also [27, 22] for more details). In the aforementioned papers, a reconstruction of the variable $\varphi$ within the cell $c_{i,j}$ is provided. In the remainder of this subsection, we apply this polynomial reconstruction procedure to our specific case of a uniform Cartesian mesh. The reconstructed variables are thus defined as follows:

$$\hat{\varphi}_{i,j}(x, y; d) = \varphi_{i,j} + \sum_{|\alpha| \in [\![1,d]\!]} R_{i,j}^\alpha \bigg( (x - x_{i,j})^{\alpha_1} (y - y_{i,j})^{\alpha_2} - M_{i,j}^\alpha \bigg), \tag{3.1}$$

where $\alpha = (\alpha_1, \alpha_2) \in \mathbb{N}^2$ is a multi-index, $|\alpha| = \alpha_1 + \alpha_2$ is its length, and $(R_{i,j}^\alpha)_{|\alpha| \in [\![1,d]\!]}$ are the polynomial coefficients.

The quantity $M_{i,j}^\alpha$ is introduced in (3.1) to ensure that the following conservation property holds for the polynomial $\hat{\varphi}_{i,j}(x, y; d)$:

$$\frac{1}{|c_{i,j}|} \int_{c_{i,j}} \hat{\varphi}_{i,j}(x, y; d) \, dx \, dy = \varphi_{i,j}.$$

Thus, we define $M_{i,j}^\alpha$ as follows:

$$M_{i,j}^\alpha = \frac{1}{|c_{i,j}|} \int_{c_{i,j}} (x - x_{i,j})^{\alpha_1} (y - y_{i,j})^{\alpha_2} \, dx \, dy.$$

Using the definition (2.32) of the cell $c_{i,j}$, we exactly compute the above integral, to get:

$$M^\alpha = \frac{1 + (-1)^{\alpha_1}}{2(\alpha_1 + 1)} \left( \frac{\Delta x}{2} \right)^{\alpha_1} \frac{1 + (-1)^{\alpha_2}}{2(\alpha_2 + 1)} \left( \frac{\Delta y}{2} \right)^{\alpha_2},$$

where we have dropped the subscript for the sake of clarity, since $M_{i,j}^\alpha$ does not actually depend on the cell $c_{i,j}$ in the uniform Cartesian situation.

Let $\Sigma_{i,j}^d$ be the stencil made of cells neighboring $c_{i,j}$ for a reconstruction of degree $d$. The stencil only depends on the degree of the polynomial reconstruction, and its construction will be detailed later on. The weights are chosen so as to minimize, in a least squares sense, the error between the average of the polynomial and the approximate solution on the cells of $\Sigma_{i,j}^d$. The polynomial coefficients $R_{i,j}^\alpha$ are therefore determined to minimize the following quadratic functional (see [21]):

$$E_{i,j}(R_{i,j}) = \sum_{l \in \Sigma_{i,j}^d} \left( \frac{1}{|c_l|} \int_{c_l} \hat{\varphi}_{i,j}(x, y; d) \, dx \, dy - \varphi_l \right)^2, \tag{3.2}$$

where $R_{i,j} = (R_{i,j}^\alpha)_{|\alpha| \in [\![1,d]\!]}$. After [21] (see also [46] for more details), we rewrite below the minimization problem (3.2) as a linear system, whose solution minimizes $E_{i,j}$.

After integrating (3.1) over the cell $c_l \in \Sigma_{i,j}^d$ and performing straightforward computations, we obtain

$$\frac{1}{|c_l|} \int_{c_l} \hat{\varphi}_{i,j}(x, y; d) \, dx \, dy = \varphi_{i,j} + \\ \sum_{|\alpha| \in [\![1,d]\!]} R_{i,j}^\alpha \left( \frac{1}{|c_l|} \left( \int_{x_l - x_{i,j} - \Delta x/2}^{x_l - x_{i,j} + \Delta x/2} \chi_{i,j}^{\alpha_1} \, d\chi_{i,j} \right) \left( \int_{y_l - y_{i,j} - \Delta y/2}^{y_l - y_{i,j} + \Delta y/2} \upsilon_{i,j}^{\alpha_2} \, d\upsilon_{i,j} \right) - M^\alpha \right), \tag{3.3}$$

where $(x_l, y_l)$ is the center of the cell $c_l$, and where we have introduced the change of variables $(\chi_{i,j}, \upsilon_{i,j}) = (x - x_{i,j}, y - y_{i,j})$. Note that there exist $(\sigma_x^l, \sigma_y^l) \in \mathbb{Z}^2$ such that $x_l - x_{i,j} = \sigma_x^l \Delta x$ and $y_l - y_{i,j} = \sigma_y^l \Delta y$. The pair of integers $(\sigma_x^l, \sigma_y^l)$ represent the position of the cell $c_l$ relatively to the position of $c_{i,j}$, as shown in Figure 4. Note that $(\sigma_x^l, \sigma_y^l)$ does not depend on the absolute position of the cell $c_{i,j}$, since our Cartesian mesh is uniform and the stencil size is the same for each cell.

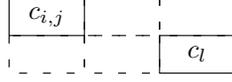

Figure 4: Relative position of the cell $c_l$ with respect to $c_{i,j}$. In this example, $\sigma_x^l = 2$ and $\sigma_y^l = -1$.

Using $\sigma_x^l$ and $\sigma_y^l$ in (3.3), we get:

$$\frac{1}{\Delta x} \int_{x_l - x_{i,j} - \Delta x/2}^{x_l - x_{i,j} + \Delta x/2} \chi_{i,j}^{\alpha_1} \, d\chi_{i,j} = \frac{(2\sigma_x^l + 1)^{\alpha_1 + 1} - (2\sigma_x^l - 1)^{\alpha_1 + 1}}{2(\alpha_1 + 1)} \left(\frac{\Delta x}{2}\right)^{\alpha_1},$$

$$\frac{1}{\Delta y} \int_{y_l - y_{i,j} - \Delta y/2}^{y_l - y_{i,j} + \Delta y/2} \upsilon_{i,j}^{\alpha_2} \, d\upsilon_{i,j} = \frac{(2\sigma_y^l + 1)^{\alpha_2 + 1} - (2\sigma_y^l - 1)^{\alpha_2 + 1}}{2(\alpha_2 + 1)} \left(\frac{\Delta y}{2}\right)^{\alpha_2}.$$

Thus, from (3.3), we obtain

$$\frac{1}{|c_l|} \int_{c_l} \hat{\varphi}_{i,j}(x, y; d) \, dx \, dy = \varphi_{i,j} +$$
$$\sum_{|\alpha| \in [\![1,d]\!]} R_{i,j}^\alpha \left( \frac{(2\sigma_x^l + 1)^{\alpha_1 + 1} - (2\sigma_x^l - 1)^{\alpha_1 + 1}}{2(\alpha_1 + 1)} \left(\frac{\Delta x}{2}\right)^{\alpha_1} \frac{(2\sigma_y^l + 1)^{\alpha_2 + 1} - (2\sigma_y^l - 1)^{\alpha_2 + 1}}{2(\alpha_2 + 1)} \left(\frac{\Delta y}{2}\right)^{\alpha_2} - M^\alpha \right). \quad (3.4)$$

Therefore, plugging (3.4) into (3.2), we have $E_{i,j}(R_{i,j}) = \|X R_{i,j} - \Phi_{i,j}\|^2$, where:

- $R_{i,j} = (R_{i,j}^\alpha)_{|\alpha| \in [\![1,d]\!]}$ is the unknown vector;
- $\Phi_{i,j} = (\varphi_l - \varphi_{i,j})_{l \in \Sigma_{i,j}^d}$;
- the matrix $X$ is defined as follows:

$$X = \left[ \frac{(2\sigma_x^l + 1)^{\alpha_1 + 1} - (2\sigma_x^l - 1)^{\alpha_1 + 1}}{2(\alpha_1 + 1)} \left(\frac{\Delta x}{2}\right)^{\alpha_1} \right.$$
$$\left. \frac{(2\sigma_y^l + 1)^{\alpha_2 + 1} - (2\sigma_y^l - 1)^{\alpha_2 + 1}}{2(\alpha_2 + 1)} \left(\frac{\Delta y}{2}\right)^{\alpha_2} - M^\alpha \right]_{l \in \Sigma_{i,j}^d, |\alpha| \in [\![1,d]\!]}.$$

Finally, we obtain $R_{i,j}$ by using the normal equation associated to the minimization problem, as follows:

$$^tXX R_{i,j} = {^tX} \Phi_{i,j}. \quad (3.5)$$

*Remark 4.* The reader is referred to [21, 26] for more details on how to efficiently solve the linear system (3.5). Indeed, since $X$ only depends on the fixed geometry, we avoid solving the linear system at each time iteration thanks to the Moore-Penrose pseudoinverse of $X$. More details on this pseudoinverse can be found in [46].

*Remark 5.* In order to ensure that there exists a solution to the minimization problem (3.2), we need more information from the stencil than we have reconstruction weights. Thus, we need $\#\Sigma_{i,j}^d > \#\{\alpha \in \mathbb{N}^2 \; ; \; |\alpha| \in [\![1,d]\!]\}$. After straightforward computations, we have the following lower bound on the size of the stencil:

$$\#\Sigma_{i,j}^d > \frac{(d+1)(d+2)}{2} - 1. \quad (3.6)$$

Therefore, to determine the stencil, we take the smallest symmetric stencil whose size satisfies (3.6) and that leads to the matrix $X$ being invertible. These choices are detailed in Figure 5, for $d \in [\![1, 5]\!]$.

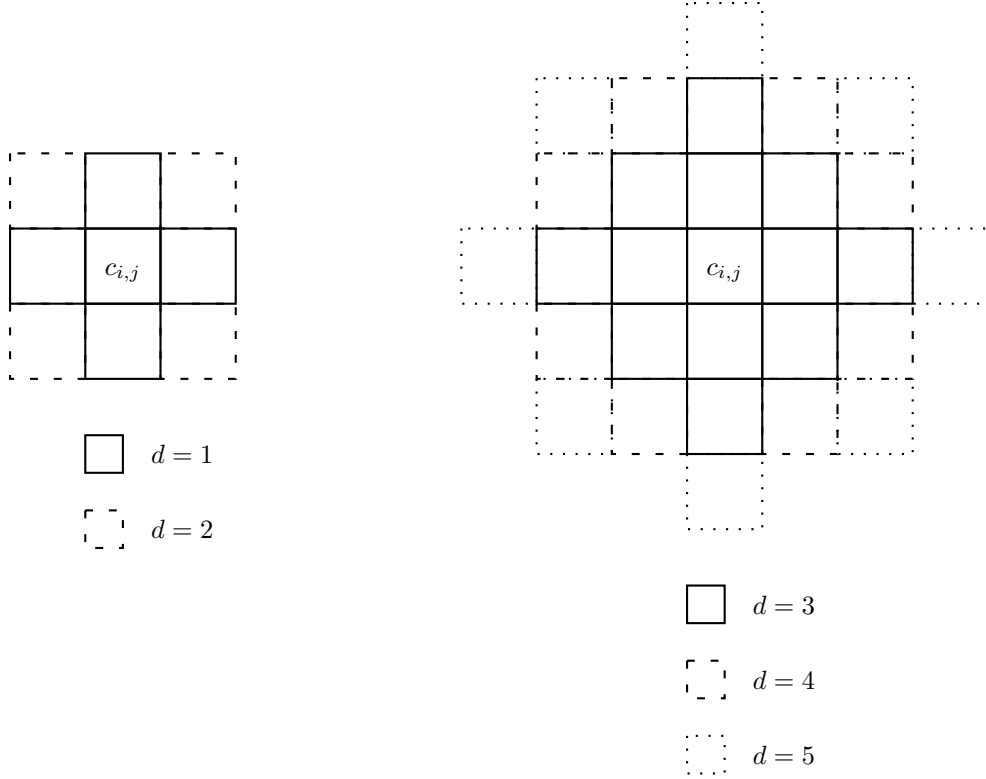

Figure 5: Representation of the stencil $\Sigma_{i,j}^d$ for $d \in [\![1,5]\!]$. The lower order stencils are always included in the higher order ones. For the sake of simplicity, we take $\Delta x = \Delta y$ in this figure.

### 3.2. The high-order scheme

Equipped with the polynomial reconstruction, out goal is now to obtain a high-order scheme, but we face a new difficulty since this high-order scheme will not be well-balanced due to the polynomial reconstruction. In a later Section, we shall introduce a correction to ensure the well-balancedness of the high-order scheme. We first present the space scheme, then its associated high-order time discretization.

#### 3.2.1. High-order space discretization

In order to improve the spatial order of accuracy of the scheme, we numerically integrate the flux at the interfaces, which requires high-order quadrature formulas based on Gauss points. The number of Gauss points $N_G$ depends only on the degree $d$, and is given by

$$N_G = 1 + \left\lfloor \frac{d}{2} \right\rfloor.$$

Let $e_{i+\frac{1}{2},j}$ be thee common interface between cells $c_{i,j}$ and $c_{i+1,j}$. The $r^{\text{th}}$ Gauss point on $e_{i+\frac{1}{2},j}$ is denoted by $\zeta_{i+\frac{1}{2},j}^r$, with the associated weight $\xi_r$ (see for instance [1] for more details on the coordinates of the Gauss points as well as their weights). Figure 6 shows the approximate location of the Gauss points on the edges of cell $c_{i,j}$ in the specific case where $N_G = 2$.

The high-order scheme reads (see [17, 26, 22] for instance):

$$
\begin{aligned}
W_{i,j}^{n+1} = W_{i,j}^n &- \sum_{r=1}^{N_G} \xi_r \left[ \frac{\Delta t}{\Delta x} \left( \mathcal{F}_{i+\frac{1}{2},j,r}^n - \mathcal{F}_{i-\frac{1}{2},j,r}^n \right) \right] \\
&- \sum_{r=1}^{N_G} \xi_r \left[ \frac{\Delta t}{\Delta y} \left( \mathcal{G}_{i,j+\frac{1}{2},r}^n - \mathcal{G}_{i,j-\frac{1}{2},r}^n \right) \right] + \Delta t (\mathcal{S}^t)_{i,j}^n + \Delta t (\mathcal{S}^f)_{i,j}^n.
\end{aligned}
\quad (3.7)
$$

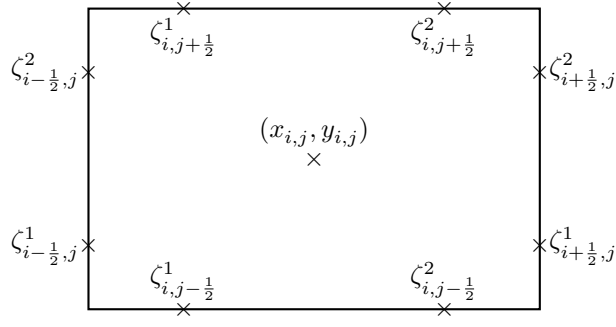

Figure 6: Approximate location of the Gauss points on the edges of cell $c_{i,j}$. We have assumed $N_G = 2$ for this figure.

The quantities $\mathcal{F}^n_{i+\frac{1}{2},j,r}$ and $\mathcal{G}^n_{i,j+\frac{1}{2},r}$ are the numerical fluxes evaluated at the edge Gauss points, given as follows:

$$\begin{aligned}
\mathcal{F}^n_{i+\frac{1}{2},j,r} &= \mathcal{F}\left(\widehat{W}^n_{i,j}(\zeta^r_{i+\frac{1}{2},j};d),\ \widehat{W}^n_{i+1,j}(\zeta^r_{i+\frac{1}{2},j};d)\right), \\
\mathcal{G}^n_{i,j+\frac{1}{2},r} &= \mathcal{G}\left(\widehat{W}^n_{i,j}(\zeta^r_{i,j+\frac{1}{2}};d),\ \widehat{W}^n_{i,j+1}(\zeta^r_{i,j+\frac{1}{2}};d)\right).
\end{aligned} \quad (3.8)$$

In (3.8), $\widehat{W}^n_{i,j}$ is the polynomial function containing the polynomial reconstructions of $h$, $q_x$ and $q_y$ within the cell $c_{i,j}$, and the functions $\mathcal{F}$ and $\mathcal{G}$ are the same as in (2.34), with the notable exception of the approximate friction source term $\overline{S}^f$ within the numerical flux, which is no longer defined by (2.20). Indeed, so as not to introduce an error in $\Delta x$ in the high-order flux approximation, we have replaced the definition (2.20) of $\overline{S}^f$ with the following expression:

$$\overline{S}^f := \overline{S}^f(h_L, h_R, q_L, q_R) = -k\bar{q}|\bar{q}|\overline{h^{-\eta}}\Delta x^d, \quad (3.9)$$

where $\bar{q}$ is defined by (2.10) and $\overline{h^{-\eta}}$ is given by:

$$\overline{h^{-\eta}} = \frac{[h^2]}{2}\frac{\eta+2}{[h^{\eta+2}]} - \frac{\bar{\mu}}{k\Delta x^{d+1}}[h]_c\left(-\frac{1}{h_L h_R} + \frac{h_L + h_R}{2}\frac{[h^{\eta-1}]}{\eta-1}\frac{\eta+2}{[h^{\eta+2}]}\right), \quad (3.10)$$

instead of (2.23). Note that, if $d = 0$, the expressions (3.9) and (3.10) coincide with (2.20) and (2.23), and the numerical flux is not modified.

In the high-order scheme (3.7), $(\mathcal{S}^t)^n_{i,j}$ and $(\mathcal{S}^f)^n_{i,j}$ are the high-order numerical source terms of topography and friction, defined as follows:

$$(\mathcal{S}^t)^n_{i,j} = \frac{1}{|c_{i,j}|}\int_{c_{i,j}}\begin{pmatrix} 0 \\ -g\hat{h}^n_{i,j}\boldsymbol{\nabla}\hat{Z}^n_{i,j}\end{pmatrix} dx\,dy, \quad (3.11a)$$

$$(\mathcal{S}^f)^n_{i,j} = \frac{1}{|c_{i,j}|}\int_{c_{i,j}}\begin{pmatrix} 0 \\ -k\hat{\boldsymbol{q}}^n_{i,j}\|\hat{\boldsymbol{q}}^n_{i,j}\|(\hat{h}^n_{i,j})^{-\eta}\end{pmatrix} dx\,dy, \quad (3.11b)$$

where $\hat{Z}^n_{i,j}$ is the reconstruction of $Z$ within the cell $c_{i,j}$, computed from the reconstructions of $h$ and $h+Z$, and where $\hat{\boldsymbol{q}}^n_{i,j} = {}^t((\hat{q}_x)^n_{i,j}, (\hat{q}_y)^n_{i,j})$. Note that we do not compute the exact integrals involved in (3.11), but rather introduce a quadrature formula of order $(d+1)$ in the cell $c_{i,j}$. The reader is referred to [1] for more information on high-order quadrature rules on a rectangle.

### 3.2.2. High-order time discretization

Strong stability-preserving Runge-Kutta (SSPRK) methods, introduced in [33, 34], are used to increase the time accuracy of the scheme, thus providing a high-order time accuracy while retaining some robustness properties of the original scheme (3.7). The second-order SSPRK(2,2), third-order SSPRK(3,3) or fourth-order SSPRK(5,4) methods, described in [33, 45], are used in the present study. Note that the SSPRK(2,2) method is nothing but Heun's method. Table 1 displays the choice of the time discretization with respect to the degree $d$.

|  $d = 1$  |  $d = 2$  |  $d \geq 3$  |
| :---: | :---: | :---: |
| SSPRK(2,2) | SSPRK(3,3) | SSPRK(5,4) |

Table 1: Choice of SSPRK method with respect to the degree of the polynomial reconstruction.

For the sake of simplicity, we shall only present the SSPRK(3,3) method, since the other two methods are similar. The reader is referred to [32] for an overview of the three SSPRK methods we use, as well as additional SSPRK methods. We begin by rewriting the scheme (3.7) as follows:

$$W^{n+1} = \mathcal{H}(W^n),$$

where $W^n$ is the vector containing the $W_{i,j}^n$ for $(i,j) \in \mathbb{Z}^2$, and $\mathcal{H}$ is a functional representing the scheme (3.7). With this notation, the SSPRK(3,3) scheme is given as follows:

$$W^{n+1} = \frac{W^n + 2W^{(3)}}{3}, \quad \text{with} \quad \begin{cases} W^{(1)} &= \mathcal{H}(W^n), \\ W^{(2)} &= \mathcal{H}(W^{(1)}), \\ W^{(3)} &= \mathcal{H}\left(\frac{3W^n + W^{(2)}}{4}\right). \end{cases} \quad (3.12)$$

The final step in the construction of the high-order time discretization is the choice of the time step $\Delta t$. For $d \leq 3$, the time step is constrained with the classical CFL condition (2.36). However, since the SSPRK(5,4) discretization is only fourth-order accurate in time, we have to introduce a correction of the time step for $d = 4$ and $d = 5$, as follows:

$$\Delta t \leq \frac{\delta^{\frac{\max(d,3)}{3}}}{2\Lambda}, \quad (3.13)$$

where $\delta$ is the 2D mesh step, given by $\delta = \min(\Delta x, \Delta y)$. The time step condition (3.13) ensures that the time scheme will be the same order of accuracy as the space scheme.

## 4. Well-balancedness recovery for the high-order scheme

The reconstruction procedure introduced in Section 3.1 causes the high-order scheme to no longer exactly preserve steady solutions, and causes non-physical oscillations when dealing with non-smooth solutions. In this section, we deal with the preservation of the steady solutions, assumed to be smooth. The oscillations are treated in the next section.

In order to recover the well-balancedness property, we propose a convex combination procedure between the first-order scheme and the high-order scheme. This convex combination recovers the well-balancedness by gradually downgrading the high-order scheme into the first-order well-balanced scheme when the solution becomes close enough to a steady state. This specific approach has been introduced in [40, 41] to produce a second-order well-balanced scheme (see also [38] for related work), and a proof of second-order accuracy was obtained in [11]. Here, the goal is to provide a very simple expression of the convex combination parameter that ensures high-order accuracy.

### 4.1. The well-balancedness property

We set up the convex combination of the first-order and the high-order schemes by introducing a pair of convex combination parameters $\boldsymbol{\theta}_{i,j}^n := {}^t((\theta_x)_{i,j}^n, (\theta_y)_{i,j}^n)$, whose expression will be given later.

Recall that the first-order well-balanced scheme is given by the two steps (2.35) – (2.38), while the high-order scheme is defined by (3.7). For the sake of clarity, we artificially split the high-order scheme (3.7) into two steps, the first one with the transport and the topography, and the second one with the friction, as follows:

$$W_{i,j}^{n+\frac{1}{2}} = W_{i,j}^n - \sum_{r=1}^{N_G} \xi_r \Delta t \left[ \frac{\mathcal{F}_{i+\frac{1}{2},j,r}^n - \mathcal{F}_{i-\frac{1}{2},j,r}^n}{\Delta x} + \frac{\mathcal{G}_{i,j+\frac{1}{2},r}^n - \mathcal{G}_{i,j-\frac{1}{2},r}^n}{\Delta y} \right] + \Delta t (\mathcal{S}^t)_{i,j}^n,$$

$$W_{i,j}^{n+1} = W_{i,j}^{n+\frac{1}{2}} + \Delta t (\mathcal{S}^f)_{i,j}^n.$$

We denote by $(\boldsymbol{S}^t)_{i,j}^n$ the second and third components of $(\mathcal{S}^t)_{i,j}^n$, defined by (3.11) and containing the topography contribution to the discharge for the high-order scheme. The first step of the high-order well-balanced scheme is merely a convex combination in each direction of the first steps of the two schemes:

$$
\begin{aligned}
W_{i,j}^{n+\frac{1}{2}} = W_{i,j}^n &- (\theta_x)_{i,j}^n \frac{\Delta t}{\Delta x} \sum_{r=1}^{N_G} \xi_r \left( \mathcal{F}_{i+\frac{1}{2},j,r}^n - \mathcal{F}_{i-\frac{1}{2},j,r}^n \right) \\
&- \left(1 - (\theta_x)_{i,j}^n\right) \frac{\Delta t}{\Delta x} \left( \mathcal{F}\left(W_{i,j}^n, W_{i+1,j}^n\right) - \mathcal{F}\left(W_{i-1,j}^n, W_{i,j}^n\right) \right) \\
&- (\theta_y)_{i,j}^n \frac{\Delta t}{\Delta y} \sum_{r=1}^{N_G} \xi_r \left( \mathcal{G}_{i,j+\frac{1}{2},r}^n - \mathcal{G}_{i,j-\frac{1}{2},r}^n \right) \\
&- \left(1 - (\theta_y)_{i,j}^n\right) \frac{\Delta t}{\Delta y} \left( \mathcal{G}\left(W_{i,j}^n, W_{i,j+1}^n\right) - \mathcal{G}\left(W_{i,j-1}^n, W_{i,j}^n\right) \right) \\
&+ \Delta t \begin{pmatrix} 0 \\ \boldsymbol{\theta}_{i,j}^n \cdot (\boldsymbol{S}^t)_{i,j}^n + \left(1 - \boldsymbol{\theta}_{i,j}^n\right) \cdot (\boldsymbol{S}^t)_{i,j}^n \end{pmatrix}.
\end{aligned}
\tag{4.1}
$$

Concerning the updated water heights, we take $h_{i,j}^{n+1} = h_{i,j}^{n+\frac{1}{2}}$, since the last step is devoted to the friction source term and therefore has no impact on the water height. We denote by $(\boldsymbol{S}^f)_{i,j}^n$ the second and third components of $(\mathcal{S}^f)_{i,j}^n$, as defined by (3.11). Following (2.38), let $(\boldsymbol{q}_{\text{WB}})_{i,j}^{n+1}$ be the vector containing the discharge obtained after the second step of the first-order scheme. The second step of the high-order well-balanced scheme consists in the convex combination of the two second steps, as follows:

$$
\boldsymbol{q}_{i,j}^{n+1} = \boldsymbol{\theta}_{i,j}^n \cdot \left( \boldsymbol{q}_{i,j}^{n+\frac{1}{2}} + \Delta t (\boldsymbol{S}^f)_{i,j}^n \right) + \left(1 - \boldsymbol{\theta}_{i,j}^n\right) \cdot (\boldsymbol{q}_{\text{WB}})_{i,j}^{n+1}.
\tag{4.2}
$$

In the two-step scheme (4.1) – (4.2), if $\boldsymbol{\theta}_{i,j}^n$ is close to 1, then the high-order scheme is favored, while the first-order well-balanced scheme is used if $\boldsymbol{\theta}_{i,j}^n$ is close to 0.

In order to recover the high-order accuracy in time, we apply the relevant SSPRK procedure with respect to the degree of the polynomial reconstruction, according to Table 1. For instance, if $d = 2$, the SSPRK(3,3) method (3.12) is used.

### 4.2. A high-order accurate convex combination

The convex combination detailed in Section 4.1 is performed in each cell $c_{i,j}$, while computing the numerical fluxes and the numerical source terms. The only remaining unknown is the pair of parameters $\boldsymbol{\theta}_{i,j}^n$. Note that, if $\boldsymbol{\theta}_{i,j}^n$ is an approximation of 1 up to order $\Delta x^{d+1}$ in each direction, then the two-step scheme with convex combination (4.1) – (4.2) is automatically of order $(d+1)$, like the fully high-order scheme (3.7).

In the remainder of this section, we derive an expression of the convex combination parameter in the $x$-direction $(\theta_x)_{i,j}^n$. The expressions in the $y$-direction are obtained in a similar fashion. Let us define a steady state detector $\varepsilon_x$ in the $x$-direction between two states $W_L$ and $W_R$, as follows:

$$
\begin{aligned}
\varepsilon_x(W_L, W_R, Z_L, Z_R, x_L, x_R) = &\sqrt{\left(\psi_x^t(W_R, Z_R) - \psi_x^t(W_L, Z_L)\right)^2 + \left((q_x)_R - (q_x)_L\right)^2 + \frac{1}{2}\left((q_y)_R^2 + (q_y)_L^2\right)} \\
&\times \sqrt{\left(\psi_x^f(W_R, x_R) - \psi_x^f(W_L, x_L)\right)^2 + \left((q_x)_R - (q_x)_L\right)^2 + \frac{1}{2}\left((q_y)_R^2 + (q_y)_L^2\right)},
\end{aligned}
\tag{4.3}
$$

where $\psi_x^t(W, Z)$ and $\psi_x^f(W, x)$ are respectively defined by the algebraic expressions of the topography and friction steady states (1.11) and (1.13), as follows:

$$
\psi_x^t(W, Z) = \frac{q^2}{2h^2} + g(h + Z),
\tag{4.4a}
$$

$$
\psi_x^f(W, x) = -q^2 \frac{h^{\eta-1}}{\eta - 1} + g \frac{h^{\eta+2}}{\eta + 2} + kq|q|x.
\tag{4.4b}
$$

Note that this definition ensures that the steady state detector $\varepsilon_x$ vanishes as soon as a steady state with either topography or friction is detected between the states $W_L$ and $W_R$ in the $x$-direction, that is to say as soon as the pairs $(h_L, (q_x)_L, Z_L)$ and $(h_R, (q_x)_R, Z_R)$ define a 1D steady state, and that $(q_y)_L = (q_y)_R = 0$.

We define the convex combination parameter in the $x$ direction $(\theta_x)_{i,j}^n$ within the cell $c_{i,j}$ as follows:

$$(\theta_x)_{i,j}^n = \sqrt{\frac{1}{2}\left(\left((\theta_x)_{i-\frac{1}{2},j}^n\right)^2 + \left((\theta_x)_{i+\frac{1}{2},j}^n\right)^2\right)}, \tag{4.5}$$

where the convex combination parameter $(\theta_x)_{i+\frac{1}{2},j}^n$ at the interface $e_{i+\frac{1}{2},j}$ is given by:

$$(\theta_x)_{i+\frac{1}{2},j}^n = \frac{(\varepsilon_x)_{i+\frac{1}{2},j}^n}{(\varepsilon_x)_{i+\frac{1}{2},j}^n + \left(\dfrac{\Delta x}{\mathcal{L}_x}\right)^k}, \tag{4.6}$$

with $(\varepsilon_x)_{i+\frac{1}{2},j}^n = \varepsilon_x(W_{i,j}^n, W_{i+1,j}^n, Z_{i,j}, Z_{i+1,j}, x_{i,j}, x_{i+1,j})$, $\mathcal{L}_x$ a characteristic length and $k \geq d+1$.

**Proposition 5.** *The convex combination parameter in the $x$ direction $(\theta_x)_{i,j}^n$, given by (4.5), satisfies the following properties:*

- $(\theta_x)_{i,j}^n$ *vanishes as soon as* $W_{i-1,j}^n$, $W_{i,j}^n$ *and* $W_{i+1,j}^n$ *define a steady state;*
- $(\theta_x)_{i,j}^n$ *is an approximation of* 1 *up to* $\Delta x^{d+1}$.

*Proof.* Equipped with the definition (4.3) of the steady state detector $\varepsilon_x$, we immediately obtain that $(\varepsilon_x)_{i+\frac{1}{2},j}^n$ vanishes as soon as $W_{i,j}^n$ and $W_{i+1,j}^n$ define a steady state. Therefore, in this case, $(\theta_x)_{i+\frac{1}{2},j}^n$ given by (4.6) also vanishes, which ensures that $(\theta_x)_{i,j}^n$ vanishes as soon as $W_{i-1,j}^n$, $W_{i,j}^n$ and $W_{i+1,j}^n$ define a steady state, and proves the first property.

The second property is proven by arguing a Taylor expansion of the expression (4.6). Indeed, note that

$$(\theta_x)_{i+\frac{1}{2},j}^n = \frac{1}{1 + \mathcal{C}_{i,j}^n \Delta x^k}, \quad \text{with} \quad \mathcal{C}_{i,j}^n = \frac{1}{(\varepsilon_x)_{i+\frac{1}{2},j}^n \mathcal{L}_x^k}.$$

Since $\mathcal{C}_{i,j}^n$ is a constant independent of $\Delta x$, we get

$$(\theta_x)_{i+\frac{1}{2},j}^n = 1 - \mathcal{C}_{i,j}^n \Delta x^k + \mathcal{O}(\Delta x^{2k}) = 1 + \mathcal{O}(\Delta x^{d+1}),$$

since $k$ is such that $k \geq d+1$. $\square$

In the $y$-direction, the convex combination parameter $(\theta_y)_{i,j}^n$ is obtained in a similar fashion. Equipped with the pair of parameters $\boldsymbol{\theta}_{i,j}^n$, the high-order well-balanced scheme (4.1) – (4.2) is complete, and its properties are summarized in the following result.

**Theorem 6.** *The scheme* (4.1) – (4.2), *equipped with the convex combination parameter defined in Section 4.2, is well-balanced by direction and high-order accurate.*

*Proof.* Since the convex combination parameter vanishes as soon as a steady state by direction is detected, the convex combination scheme (4.1) – (4.2) reverts to the first-order scheme in this case, which is well-balanced according to Theorem 4. In addition, since the convex combination parameter is nothing but an approximation of 1 up to $\Delta x^{d+1}$ according to Proposition 5, the resulting convex combination scheme is necessarily high-order accurate far from a steady state. The proof is thus concluded. $\square$

# 5. The MOOD method

The previous section proposes an adaptation the high-order scheme to recover the well-balancedness. However, this high-order accuracy comes with the loss of the non-negativity property, and the numerical solutions obtained with this scheme may present spurious oscillations around discontinuities (see [48, 39] for instance).

To address this issue, we use MOOD techniques (see [21, 26, 27] for an overview of this method, and [9, 29, 23, 22] for more recent applications, related to the shallow water equations and dry/wet transitions). Classical MOOD limiters, detailed below, are applied to preserve the non-negativity as well as to prevent the scheme from creating non-physical oscillations.

## 5.1. Overview of the MOOD method

The goal of the MOOD procedure is to recover the essential stability properties of a first-order scheme, for instance its robustness, by detecting whether the properties are verified by the high-order approximation. If this verification fails in some cell, the degree of the approximation is lowered in this cell, until the properties are satisfied. In this work, we use a more direct version of this method, by switching to the first-order scheme as soon as the verification fails, instead of progressively downgrading the polynomial degree. This choice is motivated by the fact that only the first-order scheme is well-balanced, and the well-balancedness preservation for the high-order scheme requires directly switching to the first-order scheme, as explained in Section 4.

The crux of the MOOD method lies in the choice of the properties that need to be satisfied by the high-order scheme, and in their detection. Detection criteria are commonly used within the MOOD procedure. For a more exhaustive description of these criteria, the reader is referred to [29, 22].

We introduce the notation $W^\star$ for the candidate solution, i.e. the solution obtained from $W^n$ using the high-order well-balanced scheme presented in Section 4. This candidate solution is then tested against the following criteria, to determine the cells where it is not acceptable.

### 5.1.1. Physical Admissibility Detector (PAD)

The PAD determines whether the approximate solution is out of the admissible states space $\Omega$. In the case of the shallow water equations, we check whether the water height is non-negative, and state that the PAD criterion fails within the cell $c_{i,j}$ if
$$h^\star_{i,j} < 0.$$
Let us underline that, equipped with the PAD, the high-order scheme is non-negativity preserving.

### 5.1.2. Discrete Maximum Principle (DMP)

Although the PAD ensures the non-negativity preservation, it does not prevent spurious oscillations from appearing in the vicinity of discontinuities. To address this issue, we use the DMP criterion to check for potential oscillations. Let $\nu_{i,j}$ be the set of cells connected to $c_{i,j}$ with an edge or a vertex. The DMP criterion fails if one of the following three checks is not fulfilled:

$$\begin{aligned}
\min_{l \in \nu_{i,j}} (h_l + Z_l) - \varepsilon_h &\leq h^\star_{i,j} + Z^\star_{i,j} \leq \min_{l \in \nu_{i,j}} (h_l + Z_l) + \varepsilon_h, \\
\min_{l \in \nu_{i,j}} ((q_x)_l) - \varepsilon_q &\leq (q_x)^\star_{i,j} \leq \min_{l \in \nu_{i,j}} ((q_x)_l) + \varepsilon_q, \\
\min_{l \in \nu_{i,j}} ((q_y)_l) - \varepsilon_q &\leq (q_y)^\star_{i,j} \leq \min_{l \in \nu_{i,j}} ((q_y)_l) + \varepsilon_q,
\end{aligned} \qquad (5.1)$$

where $\varepsilon_h$ and $\varepsilon_q$ are used to reduce the risk of oscillation overdetection, mainly due to floating point truncation errors. In practice, we usually take $\varepsilon_h = \varepsilon_q = \delta^3$, with $\delta = \min(\Delta x, \Delta y)$.

### 5.1.3. Detecting physical oscillations: the u2 criterion

The DMP criterion (5.1) can detect and eliminate physical extrema, thus resulting in a false positive that reduces the accuracy of the scheme. Therefore, we add another criterion to detect whether an extremum is physically admissible, namely the $u2$ criterion. It uses the second derivative of the polynomial reconstruction $\hat{\varphi}_{i,j}(x, y; 2)$. Note that, since $\hat{\varphi}_{i,j}(x, y; 2)$ is a second-degree polynomial, its second derivative is constant. In practice, we take $\varphi \in \{h + Z, q_x, q_y\}$. We then define the following curvatures on the cell $c_{i,j}$:

$$\mathcal{X}_{i,j}^{\min} = \min_{l \in \nu_{i,j}} (\partial_{xx}\hat{\varphi}_{i,j}, \partial_{xx}\hat{\varphi}_l), \qquad \mathcal{X}_{i,j}^{\max} = \max_{l \in \nu_{i,j}} (\partial_{xx}\hat{\varphi}_{i,j}, \partial_{xx}\hat{\varphi}_l),$$
$$\mathcal{Y}_{i,j}^{\min} = \min_{l \in \nu_{i,j}} (\partial_{yy}\hat{\varphi}_{i,j}, \partial_{yy}\hat{\varphi}_l), \qquad \mathcal{Y}_{i,j}^{\max} = \max_{l \in \nu_{i,j}} (\partial_{yy}\hat{\varphi}_{i,j}, \partial_{yy}\hat{\varphi}_l).$$

Equipped with the curvatures, we define three criteria, which are combined to form the $u2$ criterion (see [29, 22]). First, the plateau detector focuses on the micro-oscillations, and is defined as follows:

$$\max\left(\left|\mathcal{X}_{i,j}^{\min}\right|, \left|\mathcal{X}_{i,j}^{\max}\right|, \left|\mathcal{Y}_{i,j}^{\min}\right|, \left|\mathcal{Y}_{i,j}^{\max}\right|\right) \leq \delta. \tag{5.2}$$

Next, the local oscillation detector is given by:

$$\mathcal{X}_{i,j}^{\min}\mathcal{X}_{i,j}^{\max} \geq -\delta \quad \text{and} \quad \mathcal{Y}_{i,j}^{\min}\mathcal{Y}_{i,j}^{\max} \geq -\delta. \tag{5.3}$$

The third criterion involves a smoothness detector, given as follows, to assess whether the solution is locally smooth:

$$\frac{1}{2} \leq \frac{\min\left(\left|\mathcal{X}_{i,j}^{\min}\right|, \left|\mathcal{X}_{i,j}^{\max}\right|\right)}{\max\left(\left|\mathcal{X}_{i,j}^{\min}\right|, \left|\mathcal{X}_{i,j}^{\max}\right|\right)} \leq 1 \quad \text{and} \quad \frac{1}{2} \leq \frac{\min\left(\left|\mathcal{Y}_{i,j}^{\min}\right|, \left|\mathcal{Y}_{i,j}^{\max}\right|\right)}{\max\left(\left|\mathcal{Y}_{i,j}^{\min}\right|, \left|\mathcal{Y}_{i,j}^{\max}\right|\right)} \leq 1. \tag{5.4}$$

The $u2$ criterion is finally defined as a combination of these three detectors. Indeed, if a plateau is detected by (5.2) or if the solution is considered locally smooth by (5.4), then the DMP criterion becomes irrelevant and the $u2$ criterion succeeds. On the contrary, if a local oscillation is detected by (5.3), then the $u2$ criterion fails.

#### 5.1.4. The detector chain

Equipped with these detectors, we state the order in which the detectors are checked. To address this issue, we introduce the *Cell Polynomial Degree* (CPD) as an integer, associated to a cell $c_{i,j}$, such that $\text{CPD}(i, j) \in \{0, d\}$. If $\text{CPD}(i, j) = 0$, then the first-order scheme is used in the cell $c_{i,j}$. On the contrary, if $\text{CPD}(i, j) = d$, then the high-order scheme is used within cell $c_{i,j}$. Figure 7 displays the detector chain, and the effect of each detector on the CPD.

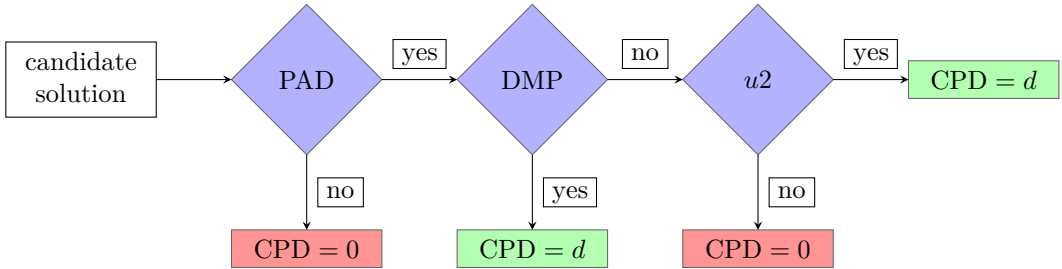

Figure 7: The MOOD detector chain.

At the end of the chain, if $\text{CPD}(i, j) = d$, the candidate solution is declared eligible in the cell $c_{i,j}$, and is accepted as the updated approximate solution $W^{n+1}$. If one of the criteria did fail, then $\text{CPD}(i, j) = 0$ is set to 0 and the candidate solution is discarded in the cell $c_{i,j}$ and its neighbors. If that is the case, a new candidate solution is computed using a polynomial reconstruction whose degree in cell $c_{i,j}$ is equal to $\text{CPD}(i, j)$. Note that, if a cell and its neighbors are declared eligible, there is no need to compute a new candidate solution in these cells.

### 5.2. Algorithm for the high-order well-balanced scheme with MOOD detection

Since the well-balancedness correction is an *a priori* procedure, it makes sense to check *a priori* for the physical admissibility of the reconstruction, in addition to using the PAD detector. The admissibility of the reconstruction is checked twice, once when computing the reconstructed heights at the Gauss points, and once when computing the numerical approximation of the mean of the friction source term. The full MOOD procedure, applied to the high-order well-balanced scheme (4.1) - (4.2), is detailed below.

**Algorithm 7.** *For a single iteration in time of the SSPRK time discretization, the MOOD loop reads as follows.*

1. For each cell $c_{i,j}$, initialize $CPD(i,j) = d$.

2. For each cell $c_{i,j}$, compute the pair of correction parameters $\boldsymbol{\theta}_{i,j}^n$. If $\boldsymbol{\theta}_{i,j}^n = \mathbf{0}$, then $CPD(i,j) = 0$.

3. For each cell $c_{i,j}$, if $CPD(i,j) > 0$, compute the interface reconstruction. If $\hat{h}_{i,j}^n(\zeta) < 0$ for some edge Gauss point $\zeta$, then the reconstruction is rejected in that cell, and we set $CPD(i,j) = 0$.

4. For each cell $c_{i,j}$, if $CPD(i,j) > 0$, compute the cell reconstruction. If $\hat{h}_{i,j}^n(\zeta) < 0$ for some cell Gauss point $\zeta$, then the reconstruction is rejected in that cell, and we set $CPD(i,j) = 0$.

5. Equipped with the new CPD map, compute the candidate solution $W^\star$, using the high-order well-balanced scheme (4.1) – (4.2).

6. Apply the detection process displayed in Figure 7 to compute a potentially new CPD map and to decide whether to accept the candidate solution. If the candidate solution is rejected, go to step 5 with the new CPD map. Otherwise, go to step 7.

7. The candidate solution is accepted, and we set $W^{n+1} = W^\star$.

## 6. Numerical experiments

This last section is devoted to numerical experiments, designed to highlight the essential properties of the scheme. The following notations are introduced to concisely label the schemes to be tested.

- The scheme that uses a polynomial reconstruction of degree $d$, i.e. whose order of accuracy is $(d+1)$, is denoted by $\mathbb{P}_d$, including the first-order well-balanced scheme.

- For $d \geq 1$, the $\mathbb{P}_d$ scheme equipped with the well-balancedness correction is denoted by $\mathbb{P}_d^{\text{WB}}$.

In addition, in order to assess the well-balancedness and the high-order accuracy of the scheme, we shall evaluate the error between the exact solution $W^{ex}(t,x,y)$ and the approximate solution. Consider a uniform Cartesian mesh made of $N = N_x \times N_y$ cells. We denote by $W_{i,j}^{ex}$ the average of the exact solution over the cell $c_{i,j}$ at time $t$, as follows:

$$W_{i,j}^{ex}(t) = \frac{1}{\Delta x \Delta y} \int_{c_{i,j}} W^{ex}(t,x,y) \, dx \, dy.$$

Equipped with this notation, we compute the errors in $L^1$, $L^2$ and $L^\infty$ norms between with $W_{i,j}^n$, the approximate solution at time $t^n$, and the exact solution $W_{i,j}^{ex}(t^n)$:

$$L^1 \text{ error: } \frac{1}{N} \sum_{i=1}^{N_x} \sum_{j=1}^{N_y} \left| W_{i,j}^n - W_{i,j}^{ex}(t^n) \right|,$$

$$L^2 \text{ error: } \sqrt{\frac{1}{N} \sum_{i=1}^{N_x} \sum_{j=1}^{N_y} \left( W_{i,j}^n - W_{i,j}^{ex}(t^n) \right)^2},$$

$$L^\infty \text{ error: } \max_{\substack{1 \leq i \leq N_x \\ 1 \leq j \leq N_y}} \left| W_{i,j}^n - W_{i,j}^{ex}(t^n) \right|.$$

The evaluation of $W_{i,j}^{ex}(t)$ for all cells $c_{i,j}$ is achieved by using a quadrature rule of the same order as the scheme (see [1] for instance). To assess the well-balancedness and the accuracy of the scheme, we evaluate these errors at the final physical time $t_{end}$.

Let us recall here that, given $\Delta x$ and $\Delta y$, the time step $\Delta t$ is given by the CFL-like condition (3.13), such that:

$$\Delta t \leq \frac{\delta^{\frac{\max(d,3)}{3}}}{2\Lambda},$$

where $\delta = \min(\Delta x, \Delta y)$ and $\Lambda$ is the maximum of all characteristic velocities at each interface.

Finally, unless otherwise specified, the two parameters from the well-balancedness detection are defined as follows: the characteristic length $\mathcal{L}_x$ is taken as the length of the domain, and we take $k = d+1$. In addition, we set $g = 9.81$ and we recall that $\eta = 7/3$.

We first propose in Section 6.1 several numerical experiments designed to assess the well-balancedness of the $\mathbb{P}_d^{\text{WB}}$ scheme, namely the preservation of the 2D lake at rest and of 1D moving steady states. Then, Section 6.2 is dedicated to the high-order accuracy on 2D exact solutions. Dam-break experiments are tackled in Section 6.3, while the simulation of the 2011 Japan tsunami is carried out in Section 6.4.

### 6.1. Well-balancedness assessment

We perform numerical experiments to assess the well-balancedness of the scheme. The first experiment concerns the preservation of a lake at rest steady state with a dry area, and the second one focuses on capturing a one-dimensional moving steady state with friction and topography that has been perturbed.

#### 6.1.1. Preservation of the lake at rest

We begin the well-balancedness numerical experiments with the preservation of a lake at rest steady state. This experiment involves a nonzero Manning coefficient $k = 10$, non-constant topography and a dry area. On the space domain $[0,1] \times [0,1]$, the topography is given by:

$$Z(x,y) = \sqrt{x^2 + y^2}.$$

The water height and the discharge are chosen according to (1.7), ensuring that $h$ stays non-negative, as follows:

$$h(t,x,y) = (1 - Z(x,y))_+ \qquad \text{and} \qquad \boldsymbol{q}(t,x,y) = \boldsymbol{0}.$$

The exact solution is prescribed as both initial and boundary conditions. A three-dimensional view of the exact height and the topography is depicted in Figure 8.

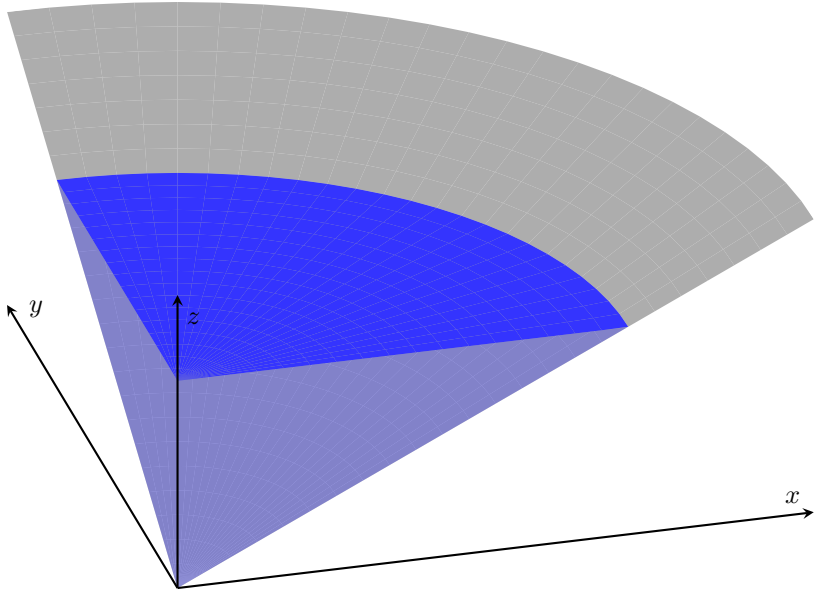

Figure 8: Exact free surface for the lake at rest experiment. The gray surface represents the topography, the opaque blue surface is the water surface and the transparent blue volume is filled with water.

In order to highlight the relevance of the well-balancedness correction, the simulation is carried out using the first-order scheme and the sixth-order scheme, with and without correction. The results of the experiment are reported in Table 2, for 2500 ($50 \times 50$) cells and at time $t_{end} = 0.1$s. For this simulation, we set the cutoff constant $C$ defined in (2.17) to $C = +\infty$. Numerically, we take $C$ as the upper bound of the double precision floating point numbers.

|         | $h + Z$    |          |          | $\|\boldsymbol{q}\|$ |          |          |
|---------|------------|----------|----------|----------|----------|----------|
|         | $L^1$      | $L^2$    | $L^\infty$ | $L^1$    | $L^2$    | $L^\infty$ |
| $\mathbb{P}_0$ | 5.76e-18 | 2.53e-17 | 2.22e-16 | 5.68e-17 | 1.05e-16 | 7.44e-16 |
| $\mathbb{P}_5$ | 5.85e-05 | 1.60e-04 | 1.51e-03 | 7.98e-05 | 2.62e-04 | 7.54e-03 |
| $\mathbb{P}_5^{\text{WB}}$ | 2.71e-17 | 6.95e-17 | 6.66e-16 | 2.59e-16 | 3.77e-16 | 2.13e-15 |

Table 2: Free surface and discharge norm errors for the lake at rest experiment.

From Table 2, we observe that the first-order scheme indeed exactly preserves the lake at rest. However, the sixth-order $\mathbb{P}_5$ scheme, as expected, does not exactly preserve the lake at rest but instead gives a sixth-order approximation of this steady state. The relevance of the correction is thus highlighted here, since it allows to recover the exact lake at rest steady state.

### 6.1.2. Subcritical steady flow over a bump without friction

To assess the well-balancedness by direction, we consider the well-known Goutal and Maurel test case from [35] in the $x$-direction. The initial conditions for this experiment consist in a lake at rest over a topography with a bump, given on the domain $(0, 25) \times (0, 1)$ by:

$$Z(x, y) = \max\left(\left[0.2 - 0.05(x - 10)^2\right], 0\right).$$

The initial free surface is given by $h(0, x, y) + Z(x, y) = 2$, and the initial discharge is set to zero. The friction is canceled, and we take $k = 0$.

The main feature of this experiment is that the final steady state is obtained following a transient state, which is governed by the boundary conditions. Neumann boundary conditions are prescribed on each boundary and each variable, except on the left boundary where the $x$-discharge is such that $q_x(t, 0, y) = 4.42$, and on the right boundary where the height is set to $h(t, 25, y) = 2$. These boundary conditions enable the eventual formation of a subcritical moving steady state. Therefore, this test case not only checks whether the well-balanced scheme under consideration is able to preserve a given steady state, but also if it is possible to capture a steady state obtained after a transient state. In [40], the first-order well-balanced scheme was shown to capture this subcritical steady state, and gives evidence that the well-balancedness correction of the $\mathbb{P}_5^{\text{WB}}$ scheme should also capture the steady solution.

The simulation is carried out on $100 = 100 \times 1$ cells, using the first-order scheme and the sixth-order scheme, with and without correction. Regarding the correction, we take $\mathcal{L}_x = 0.5$. The errors are presented in Table 3 at the final time $t_{end} = 500\text{s}$, where $\psi_x^t$ is the topography steady state detector in the $x$-direction given by (4.4a), which becomes constant once a steady state is reached. We correctly recover the expected behavior, that is to say both the $\mathbb{P}_0$ and $\mathbb{P}_5^{\text{WB}}$ schemes capture the subcritical steady state up to the machine precision, while the $\mathbb{P}_5$ scheme merely approximates this steady state. Of course, the same conclusion is reached by considering the experiment in the $y$-direction.

|         | $\psi_x^t$ |          |          | $\|\boldsymbol{q}\|$ |          |          |
|---------|------------|----------|----------|----------|----------|----------|
|         | $L^1$      | $L^2$    | $L^\infty$ | $L^1$    | $L^2$    | $L^\infty$ |
| $\mathbb{P}_0$ | 1.03e-13 | 1.85e-13 | 1.41e-13 | 6.69e-14 | 7.74e-14 | 1.41e-13 |
| $\mathbb{P}_5$ | 2.14e-04 | 1.03e-03 | 9.92e-03 | 2.12e-04 | 5.44e-04 | 2.11e-03 |
| $\mathbb{P}_5^{\text{WB}}$ | 1.05e-12 | 1.88e-12 | 4.34e-12 | 6.75e-13 | 8.25e-13 | 1.53e-12 |

Table 3: Errors on $\psi_x^t$ and on $\|\boldsymbol{q}\|$ for the subcritical Goutal and Maurel test case, obtained after a transient state.

### 6.1.3. Perturbed steady state with friction over a flat topography

We now consider a numerical experiment designed to check the well-balancedness by direction with respect to the friction source term. To that end, we take $k = 1$ and a flat topography $Z(x, y) = 0$. The initial water height derives from solving the equation $\psi_x^f = 0.02$, with $q_x(x, y) = -0.5$ and $q_y(x, y) = 0$, where $\psi_x^f$, defined by (4.4b), is constant for a friction steady state. This initial condition $W^{\text{steady}}(x, y)$, therefore, represents a steady state at rest for the friction source term, and it should be preserved by the well-balanced schemes.

To study the capture, rather than the preservation, of a friction steady state we add a perturbation to the steady solution $W^{\text{steady}}(x, y)$. The height is perturbed as follows:

$$\begin{cases} h(0, x, y) = h^{\text{steady}}(x, y) + 0.05 & \text{if } x \in \left(\frac{3}{7}, \frac{4}{7}\right), \\ h(0, x, y) = h^{\text{steady}}(x, y) & \text{otherwise,} \end{cases}$$

while the following perturbation is applied to the $x$-discharge:

$$\begin{cases} q_x(0, x, y) = q_x^{\text{steady}}(x, y) + 0.5 & \text{if } x \in \left(\frac{3}{7}, \frac{4}{7}\right), \\ q_x(0, x, y) = q_x^{\text{steady}}(x, y) & \text{otherwise.} \end{cases}$$

Equipped with Dirichlet boundary conditions given by the steady solution, the perturbation will therefore eventually be dissipated. Thus, the test case enables to test the capture of the resulting friction-only steady state, obtained after a transient state. A similar experiment was performed in [41], where the first-order scheme was shown to capture this steady state at machine precision.

We carry out the simulation on $100 = 100 \times 1$ cells with the $\mathbb{P}_0$, $\mathbb{P}_5$ and $\mathbb{P}_5^{\text{WB}}$ schemes, until the final time $t_{end} = 5$s. For the $\mathbb{P}_5^{\text{WB}}$ scheme, we take $\mathcal{L}_x = 1/15$. As expected, the $\mathbb{P}_0$ and $\mathbb{P}_5^{\text{WB}}$ scheme capture the steady solution up to machine precision. The $\mathbb{P}_5$ scheme gives an approximation of the solution, but remains far from the exact solution. Same conclusions are obtained for the $y$-direction.

|  | $\psi_x^f$ | | | $\|q\|$ | | |
| --- | --- | --- | --- | --- | --- | --- |
|  | $L^1$ | $L^2$ | $L^\infty$ | $L^1$ | $L^2$ | $L^\infty$ |
| $\mathbb{P}_0$ | 9.18e-15 | 1.08e-14 | 1.96e-14 | 3.43e-15 | 3.95e-15 | 6.88e-15 |
| $\mathbb{P}_5$ | 1.55e-03 | 1.55e-03 | 1.59e-03 | 1.04e-03 | 1.04e-03 | 1.06e-03 |
| $\mathbb{P}_5^{\text{WB}}$ | 6.53e-13 | 8.38e-13 | 1.44e-12 | 3.38e-13 | 3.82e-13 | 5.47e-13 |

Table 4: Errors on $\psi_x^f$ and on $\|q\|$ for the steady state with friction, obtained after a transient state resulting from a perturbation.

In addition, we display the results of the $\mathbb{P}_0$ and $\mathbb{P}_5^{\text{WB}}$ schemes in Figure 9, for several physical times ($t = 0$s, $t = 0.05$s and $t = 5$s). We observe that, starting with the same initial condition, we end up with the same steady state up to machine precision. However, the approximation of the transient state in the middle panel is much less diffusive when using the $\mathbb{P}_5^{\text{WB}}$ scheme compared to the $\mathbb{P}_0$ scheme. This highlight the relevance of both the high-order accuracy and the well-balancedness property.

### 6.2. Order of accuracy assessment

We now assess the accuracy of the high-order scheme. To that end, we propose two numerical experiments. The first one only involves the topography source term, while both topography and friction are considered in the second benchmark.

#### 6.2.1. Steady vortex

This first experiment is a steady vortex (see [40, 22]). On the space domain $(-1, 1)^2$, we set the topography as $Z(x, y) = 0.2 e^{0.5(1-r^2)}$, with $r^2 = x^2 + y^2$. We cancel the friction term by taking $k = 0$. The exact solution, displayed in Figure 10, is then given by $W_{ex} = {}^t(h, hu, hv)$, where we have set

$$h(t, x, y) = 1 - \frac{1}{4g} e^{2(1-r^2)} - Z(x, y) \; ; \quad u(t, x, y) = y \, e^{1-r^2} \; ; \quad v(t, x, y) = -x \, e^{1-r^2}.$$

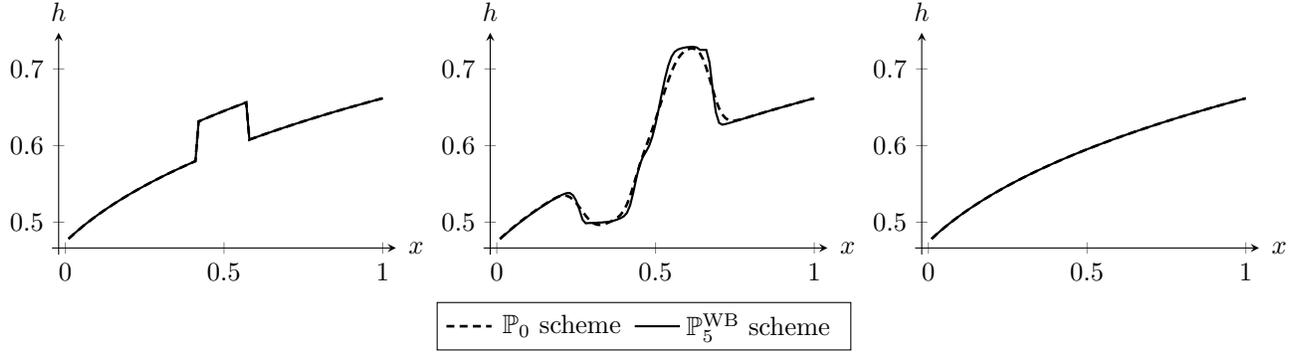

Figure 9: Water height for the perturbed steady state with friction, with the $\mathbb{P}_0$ scheme (dashed line) and the $\mathbb{P}_5^{\text{WB}}$ scheme (solid line). Left panel: initial condition at $t = 0$s; middle panel: transient state at $t = 0.05$s; right panel: steady state at $t = 5$s.

The initial conditions consist in computing a numerical average with a quadrature formula of order $(d+1)$, of the exact solution in each cell. Similarly, the boundary conditions are obtained by evaluating the exact solution at the Gauss points on the domain boundary.

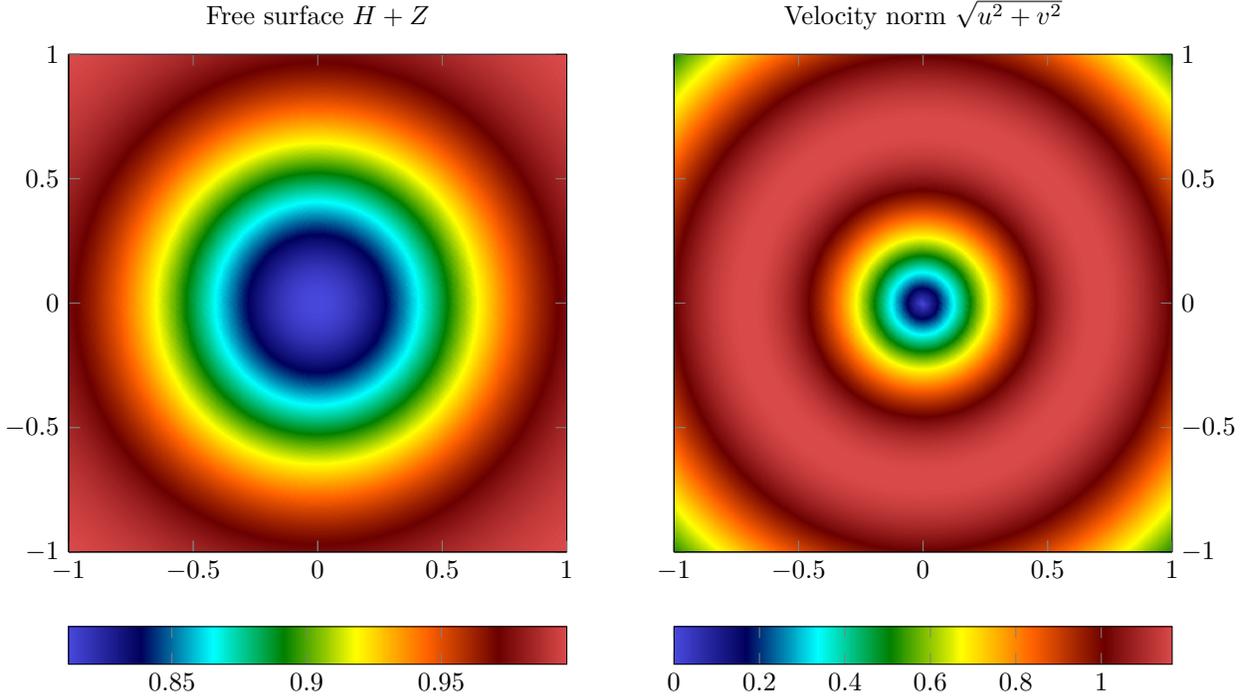

Figure 10: Exact solution for the steady vortex experiment. Left panel: free surface. Right panel: velocity norm (the vortex flows clockwise).

The simulations are carried out with the $\mathbb{P}_3^{\text{WB}}$ and $\mathbb{P}_5^{\text{WB}}$ schemes, until a final physical time $t_{end} = 1$s. In addition, we take $C = +\infty$ for each scheme. The results of the simulations are presented in Figure 11.

In all cases, the accuracy reaches the order $(d+1)$, as expected. Optimal accuracy is maintained thanks to the $u2$ detection criteria. Indeed, on such smooth solutions, the DMP criterion (5.1) would wrongly lower the CPD in some cells by over-detecting smooth extrema. Here, the smoothness detector (5.4) is used to correct over-detection from the DMP criterion. The reader is referred to [22] for a comparison of the order with and without the $u2$ criterion. In [22], the authors show that the $u2$ criterion is mandatory to recover the expected

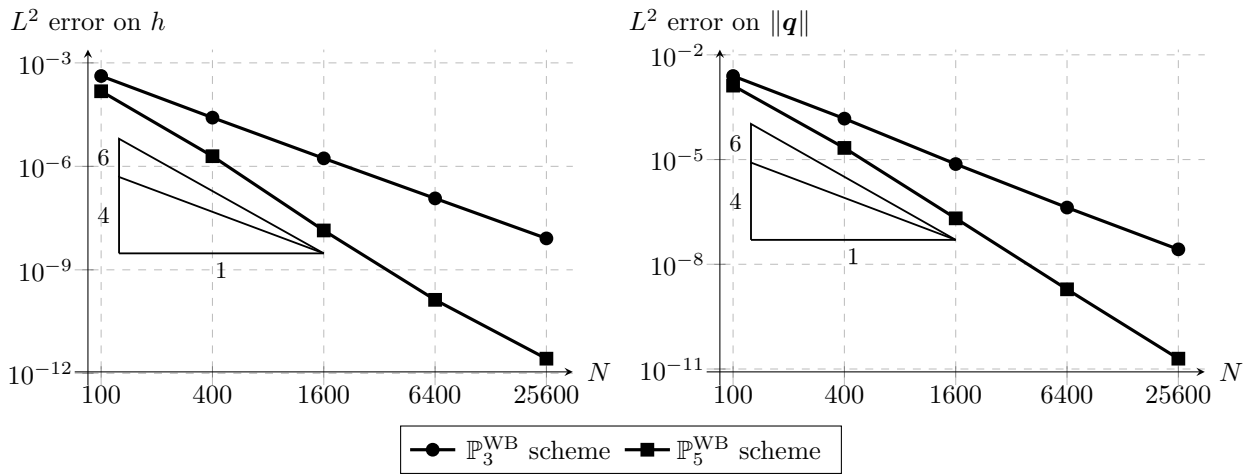

Figure 11: Error lines in $L^2$ norm for the exact solution with topography, using the $\mathbb{P}_3^{\mathrm{WB}}$ and $\mathbb{P}_5^{\mathrm{WB}}$ schemes. Left panel: $L^2$ error on $h$; right panel: $L^2$ error on $\|\boldsymbol{q}\|$.

order of accuracy.

#### 6.2.2. Experiment with topography and friction

The second experiment for accuracy assessment concerns a two-dimensional steady state involving both the topography and the friction source term (see [41]). The exact solution is defined as follows. Let $\boldsymbol{r} = {}^t(x, y)$. We assume, for this experiment, that $\|\boldsymbol{r}\| \neq 0$. The topography is given by:

$$Z(x, y) = \frac{2k\|\boldsymbol{r}\| - 1}{2g\|\boldsymbol{r}\|^2}.$$

In addition, we set $W_{ex} = {}^t(h, \boldsymbol{q})$, where

$$h(t, x, y) = 1 \qquad \text{and} \qquad \boldsymbol{q}(t, x, y) = \frac{\boldsymbol{r}}{\|\boldsymbol{r}\|^2}.$$

For the purpose of the simulation, we consider the exact solution on the space domain $(0.4, 1)^2$, with a Manning coefficient $k = 1$. As in the previous experiment, the initial and boundary conditions derive from the exact solution.

In order to check the high-order accuracy of the schemes, the benchmark is carried out with the $\mathbb{P}_3^{\mathrm{WB}}$ and $\mathbb{P}_5^{\mathrm{WB}}$ schemes. The final physical time is $t_{end} = 0.1$s, and we take once again $C = +\infty$. Convergence curves are presented in Figure 12.

Once again, we optimal the expected order of accuracy around $d + 1$. Similarly to the previous experiment, this order of accuracy is recovered thanks to the $u2$ criterion in addition to the DMP criterion.

#### 6.3. Dam-break test cases

This section is dedicated to the validation of dam-break benchmarks. We first consider a dam-break over a dry bottom in one space direction. Such a simulation will highlight the relevance of the well-balancedness correction and the MOOD procedure. Next, we present a two-dimensional partial dam-break.

#### 6.3.1. One-dimensional dry dam-break

We consider the academic square domain $[0, 1]^2$, and the topography is given by:

$$Z(x, y) = \frac{e^x}{e^1} - e^{-1},$$

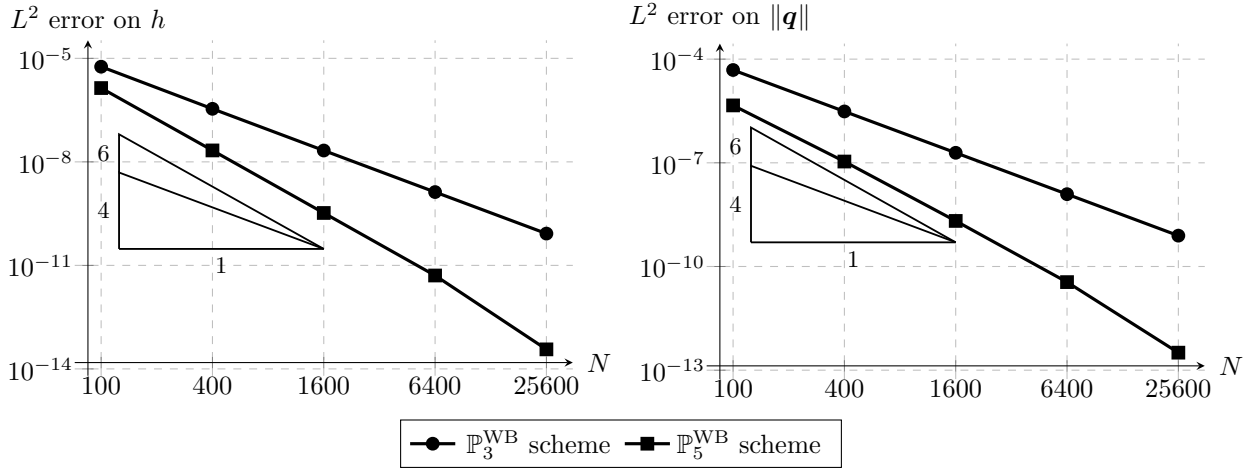

Figure 12: Error lines in $L^2$ norm for the exact solution with topography and friction, using the $\mathbb{P}_3^{\text{WB}}$ and $\mathbb{P}_5^{\text{WB}}$ schemes. Left panel: $L^2$ error on $h$; right panel: $L^2$ error on $\|\boldsymbol{q}\|$.

such that $Z(0, y) = 0$ and $Z(1, y) = 1$. The initial free surface consists in a double dam-break, obtained by setting:

$$h(0, x, y) + Z(x, y) = \begin{cases} 2 & \text{if } x < \frac{1}{2}, \\ Z(x, y) & \text{otherwise.} \end{cases}$$

In addition, the initial discharge is zero, i.e. $\boldsymbol{q}(0, x, y) = \boldsymbol{0}$. The Manning coefficient is set to 1 and the boundaries are considered to be solid walls. The experiment is carried out with the $\mathbb{P}_0$ and $\mathbb{P}_5^{\text{WB}}$ schemes. The final physical time is $t_{end} = 0.07$s, and we set $C = +\infty$. We also take $\mathcal{L}_x = 0.1$. The results are presented in Figures 13 and 14.

Figure 13 shows a comparison between the free surface obtained with the $\mathbb{P}_0$ scheme and the one obtained with the $\mathbb{P}_5^{\text{WB}}$ scheme, using $50 = 25 \times 2$ cells in each case. We also display a reference solution, obtained by using the $\mathbb{P}_0$ scheme with $1000 = 500 \times 2$ discretization cells.

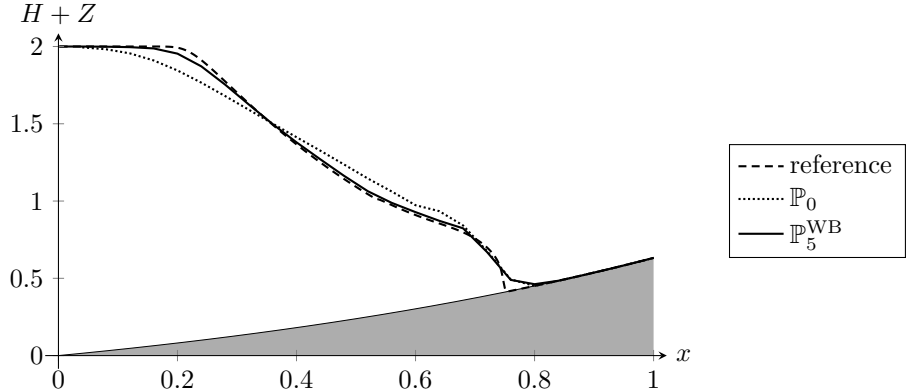

Figure 13: Free surface for the dam-break on a dry slope experiment; reference solution (dashed line), $\mathbb{P}_0$ scheme (dotted line), $\mathbb{P}_5^{\text{WB}}$ scheme (solid line). The gray area represents the topography.

Figure 13 highlights the relevance of the well-balancedness property as well as the high-order accuracy. First, despite the coarse grid, the results from the $\mathbb{P}_5^{\text{WB}}$ scheme are close to the reference solution, except in one cell close to the dry/wet transition, where the PAD detector is activated. In addition, we note that the free surface is unperturbed close to the left edge of the domain. Indeed, the waves from the dam-break have not yet reached the edges of the domain at $t = t_{end}$, and the area located in the vicinity of the left edge is in a lake at rest configuration. This essential property exactly holds for the $\mathbb{P}_5^{WB}$ scheme. This behavior is obtained thanks to the well-balancedness correction, which forces the well-balanced scheme to be activated in lake at rest-type situations.

In Figure 14, the well-balancedness coefficient in the $x$-direction $\theta_x$ is depicted, together with the free surface and the topography, for $t = t_{end}/2$ and $t = t_{end}$, for the $\mathbb{P}_5^{WB}$ scheme.

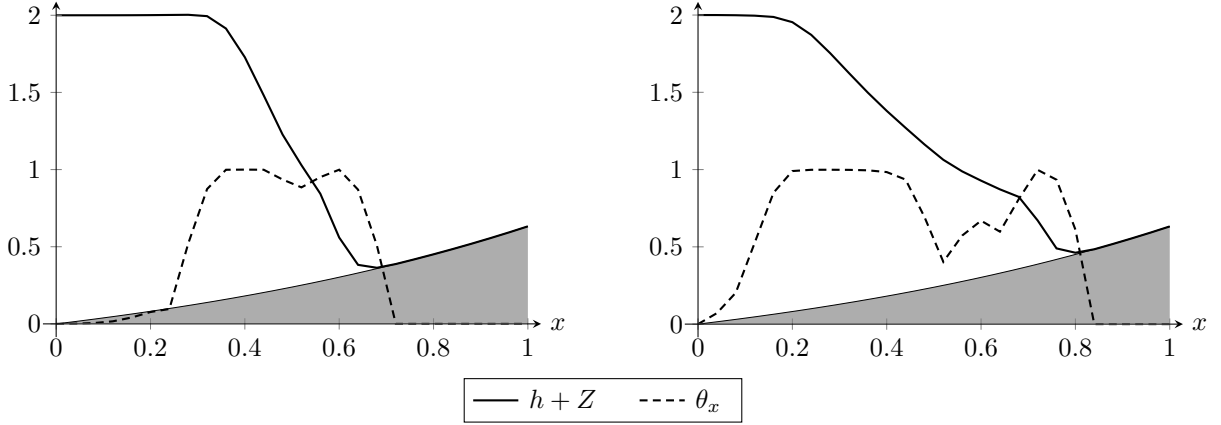

Figure 14: Free surface $h + Z$ (solid line) and well-balancedness parameter in the $x$-direction $\theta_x$ (dashed line) for the dam-break on a dry slope experiment using the $\mathbb{P}_5^{WB}$ scheme. The gray area represents the topography. Left panel: $t = t_{end}/2$; right panel: $t = t_{end}$.

In the left panel of Figure 14, we observe that $\theta_x$ is zero in areas that have not yet been impacted by the waves, i.e. in the areas where a lake at rest configuration is found. As a consequence, in these areas (close to the edges of the domain), the well-balanced scheme is used. Similar conclusions are drawn from the right panel of Figure 14. The edges of the domain are still considered to be at rest, which is evidenced by the convex combination parameter being very close to zero.

*6.3.2. Two-dimensional partial dam-break*

The second experiment concerns a two-dimensional partial dam-break (see [42, 22, 41]). An extensive study of this experiment, focusing on the differences between various reconstruction degrees and MOOD criteria, has been carried out in [22]. In [22], the authors show that the depth of the vortices appearing at the edges of the dam strongly depend on the degree of the reconstruction and the MOOD criteria used. However, in [22], the friction source term was not present, and the authors only studied the effects of the topography. Thus, in the present paper, we focus on the impact of the friction source term, by carrying out the simulation with three different Manning coefficients.

For this experiment, the space domain is $[-100, 100] \times [-100, 100]$, and the topography is given as follows:

$$Z(x, y) = \begin{cases} 1 & \text{if } x \leq -5, \\ 0 & \text{if } x \geq 5, \\ 0.1(5 - x) & \text{if } -5 < x < 5 \text{ and } -40 < y < 40, \\ 12 & \text{if } -5 < x < 5 \text{ and } y \in [-100, -40] \cup [40, 100]. \end{cases}$$

It represents a 12 meters high, 10 meters wide broken dam. Initially, the reservoir (to the left) is filled, as follows:

$$h(0, x, y) = \begin{cases} 10 - Z(x, y) & \text{if } x \leq -5, \\ 5 - Z(x, y) & \text{if } x \geq 5, \\ 5 - Z(x, y) & \text{if } -5 < x < 5 \text{ and } -40 < y < 40, \\ 0 & \text{if } -5 < x < 5 \text{ and } y \in [-100, -40] \cup [40, 100]. \end{cases}$$

The water is initially at rest, i.e. $\boldsymbol{q}(0, x, y) = \boldsymbol{0}$. For this simulation, we use wall boundary conditions. All the simulations are carried out with $40000 = 200 \times 200$ discretization cells.

The goal of this simulation is to compare the results from the $\mathbb{P}_0$, $\mathbb{P}_1^{\text{WB}}$ and $\mathbb{P}_5^{\text{WB}}$ schemes. Moreover, the simulation is carried out with various Manning coefficients, namely $k = 0$, $k = 0.25$ and $k = 2$, and until the final physical time $t_{\text{end}} = 7$s. In addition, we set $C = 0.5$.

The results of the simulations are displayed in Figure 15 ($k = 0$), in Figure 16 ($k = 0.25$) and in Figure 17 ($k = 2$). We use the same color scale in all figures.

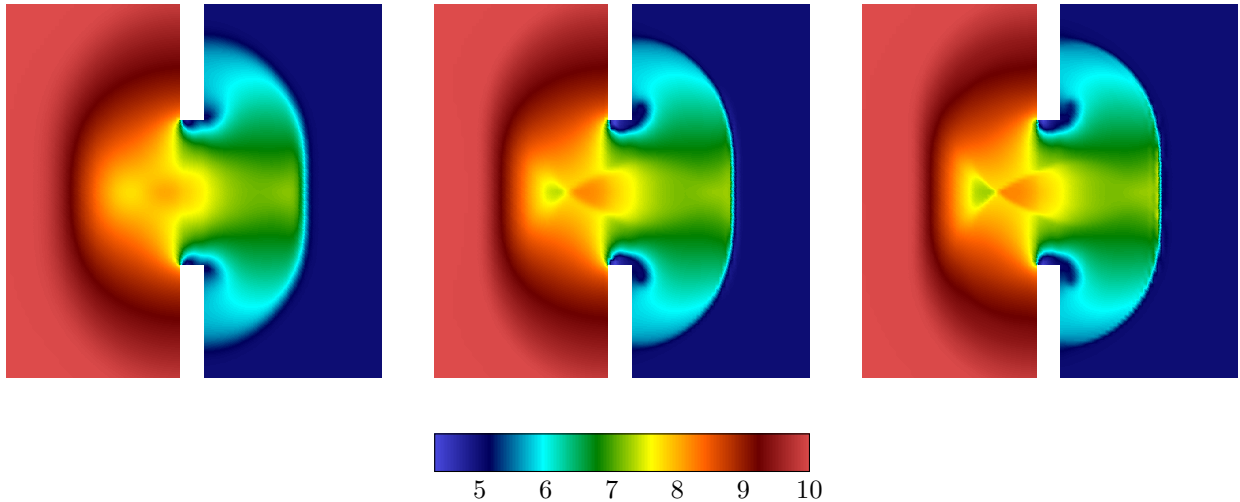

Figure 15: Free surface for the partial dam-break simulation with $k = 0$. From left to right: results of the $\mathbb{P}_0$, $\mathbb{P}_1^{\text{WB}}$ and $\mathbb{P}_5^{\text{WB}}$ schemes.

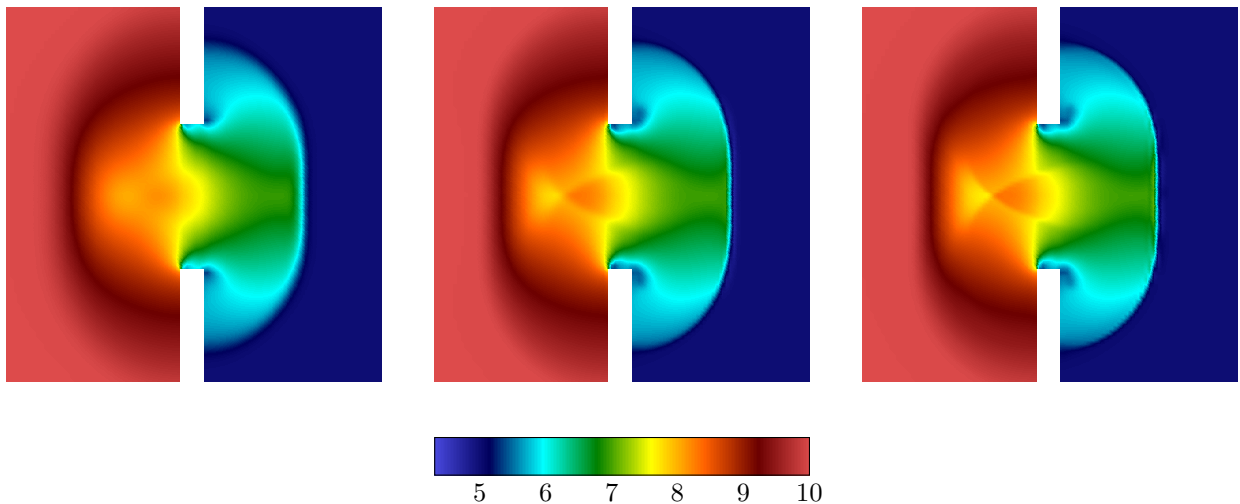

Figure 16: Free surface for the partial dam-break simulation with $k = 0.25$. From left to right: results of the $\mathbb{P}_0$, $\mathbb{P}_1^{\text{WB}}$ and $\mathbb{P}_5^{\text{WB}}$ schemes.

Figure 15 and Figure 16 show that the shock wave to the right of the dam and the rarefaction wave to the left of the dam are clearly more smeared when using the $\mathbb{P}_0$ scheme instead of the $\mathbb{P}_1^{\text{WB}}$ or the $\mathbb{P}_5^{\text{WB}}$ scheme. In addition, the shock structure at the center of the water flow is not visible with the $\mathbb{P}_0$ scheme. This structure, although smeared, is visible with the $\mathbb{P}_1^{\text{WB}}$ scheme, and turns out to be very well-defined with the $\mathbb{P}_5^{\text{WB}}$ scheme. We draw similar conclusions from Figure 17. The smearing of the shock wave and the rarefaction wave is noticeable with the first- and second-order schemes, but it is strongly reduced with the high-order scheme. In addition, the important friction has caused the central structure to nearly disappear.

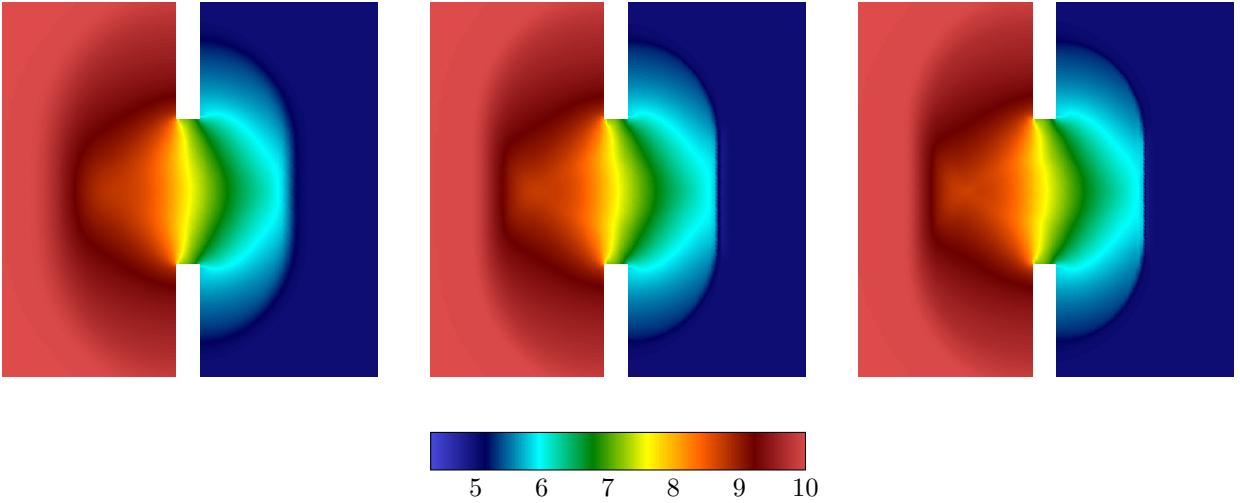

Figure 17: Free surface for the partial dam-break simulation with $k = 2$. From left to right: results of the $\mathbb{P}_0$, $\mathbb{P}_1^{\mathrm{WB}}$ and $\mathbb{P}_5^{\mathrm{WB}}$ schemes.

An important remark we make here concerns the vortices present at the edges of the dam in Figures 15 and 16. The presence of the friction source term dampens the depth, as well as the size, of these vortices. We focus on the top vortex, whose characteristics are similar to the bottom one since the experiment is symmetric with respect to the $y = 0$ line. This behavior is displayed in Table 5, where the approximate size and the depth of the vortex are collected.

| Manning coefficient | Vortex size | Vortex depth |
|---|---|---|
| $k = 0$ | $84\mathrm{m}^2$ | 4.28m |
| $k = 0.25$ | $17\mathrm{m}^2$ | 5.45m |
| $k = 2$ | $0\mathrm{m}^2$ | 7.23m |

Table 5: Depth and approximate size of the deepest vortex, for the $\mathbb{P}_5^{\mathrm{WB}}$ scheme. For the case where $k = 2$, there is no vortex, and the table displays the free surface at the point where the vortex would be located if the Manning coefficient were lower.

Concerning the left rarefaction wave, the relevant indicators are the position of the head of the rarefaction wave, its size, and its amplitude along the $y = 0$ line. Those quantities are reported in Table 6, where we chose to compute the amplitude of the rarefaction wave by subtracting the water height at the tail from the water height at the head.

| Manning coefficient | Size | Amplitude | Head |
|---|---|---|---|
| $k = 0$ | 39m | 2.68m | $x = -74$m |
| $k = 0.25$ | 38m | 2.28m | $x = -74$m |
| $k = 2$ | 31m | 1.29m | $x = -74$m |

Table 6: Left rarefaction wave: approximate size, water height amplitude and position of the head, with respect to the Manning coefficient.

Concerning the shock wave, we report on its position and its amplitude along the line $y = 0$, given in Table 7. Similarly to the rarefaction wave, the amplitude of the shock wave is obtained by computing the difference between the water height to the left of the wave and the water height to its right. Note that, since those computations are performed on the numerical results of the $\mathbb{P}_5^{\mathrm{WB}}$ scheme, the shock wave takes only a couple of cells, and the evaluation of its position is fairly accurate. In addition, the amplitude of the shock wave presented for $k = 0$ in Table 7 is very similar to the results obtained in [22], although the authors do not use the same scheme.

| Manning coefficient | Position | Amplitude |
|---|---|---|
| $k = 0$ | $x = 60$m | 2.28m |
| $k = 0.25$ | $x = 58$m | 1.96m |
| $k = 2$ | $x = 53$m | 0.98m |

Table 7: Right shock wave: approximate position and water height amplitude, with respect to the Manning coefficient.

Tables 6 and 7 give evidence about the effect of the friction on the water flow. The Manning term dampens the amplitudes of both the rarefaction wave and the shock wave, while an increase in the friction coefficient is accompanied by a diminution of the size of the rarefaction wave, and a decrease in the distance traveled by the shock wave. This behavior is expected, as an increase in friction leads to a decrease in discharge, as evidenced by the expressions (2.29) in 1D and (2.37) in 2D. The discharge decrease leads to a slower travel time of the shock wave, which means that the wave will travel less distance.

Finally, we observe from Table 6 that the friction does not change the position of the head of the rarefaction wave. This behavior is also expected from the expression of the friction source term given by (1.8) in 1D and (1.1) in 2D. Near the head of the rarefaction wave, the water is almost at rest, since no wave has already perturbed the initial rest condition, leading to a negligible impact of the friction source term, which means the head of the rarefaction wave travels at the same speed for $k = 0$, $k = 0.25$ or $k = 2$. Therefore, the value of the Manning coefficient does not alter the position of the head of the rarefaction wave.

### 6.4. Simulation of the 2011 Japan tsunami

We tackle the simulation of the 2011 Japan tsunami, based on real data, see for instance [23]. This real data consists in a uniform Cartesian mesh made of about $13 \times 10^6$ rectangles, where the cell topography and the initial free surface for the tsunami simulation are given, see Figure 18. The initial discharge is set to zero, that is $\boldsymbol{q}(0, x, y) = \boldsymbol{0}$. Homogeneous Neumann boundary conditions are prescribed at each boundary. In addition, we set $k = 0.05$.

To assess the simulation accuracy, we compare real physical measurements from sensors far from the Japanese coast with the numerical approximation. The positions of these sensors are depicted in the left panel of Figure 19. The sensors have measured the water height during one hour, and thus we set the final physical time to $t_{end} = 3600$s. The main difficulty of this simulation lies in the large topography gradients in the mesh with respect to the characteristic cell size. For instance, the right panel of Figure 19 shows the topography over the solid horizontal line in the left panel of Figure 19. The extreme topography gradients, especially around the Kuril trench, have to be correctly handled by the scheme.

We check the numerical results of the $\mathbb{P}_0$, $\mathbb{P}_1^{WB}$ and $\mathbb{P}_3^{WB}$ schemes. We ran the simulation until the final time $t_{end} = 1$h on 48 computational cores; the $\mathbb{P}_0$ scheme took around 1 hour of CPU time, whereas the $\mathbb{P}_1^{WB}$ and $\mathbb{P}_3^{WB}$ schemes took respectively around 2.5h and around 10h. The numerical results are displayed in Figure 20. The $\mathbb{P}_0$ and $\mathbb{P}_1^{WB}$ schemes yield comparable results, and the second-order result is, as expected, less diffusive than the first-order one, with much more structure present within the waves. Unfortunately, the results of the $\mathbb{P}_3^{WB}$ scheme are unsatisfactory. The extreme topography gradients present in the domain, like the ones depicted in the right panel of Figure 19, have led to an over-limitation of the MOOD method to remove the spurious oscillations, and the fourth-order solution ends up looking very similar to the first-order one.

In Figure 21, we display the sea surface height (SSH), that is to say the difference between the water height and the average surface elevation, at each of the three sensors. The physical data is compared to the results from the $\mathbb{P}_0$ and $\mathbb{P}_1^{WB}$ schemes. We observe that, although the $\mathbb{P}_0$ scheme already gives a good approximation of the data, the $\mathbb{P}_1^{WB}$ approximation is even better. Namely, the correct tsunami propagation time is captured, and the well-balancedness of the schemes ensure that no spurious oscillations come from the balance between flux and topography. This result questions the need to even use higher-order schemes for this simulation with such large cells, since the second-order results are already very close to the physical data.

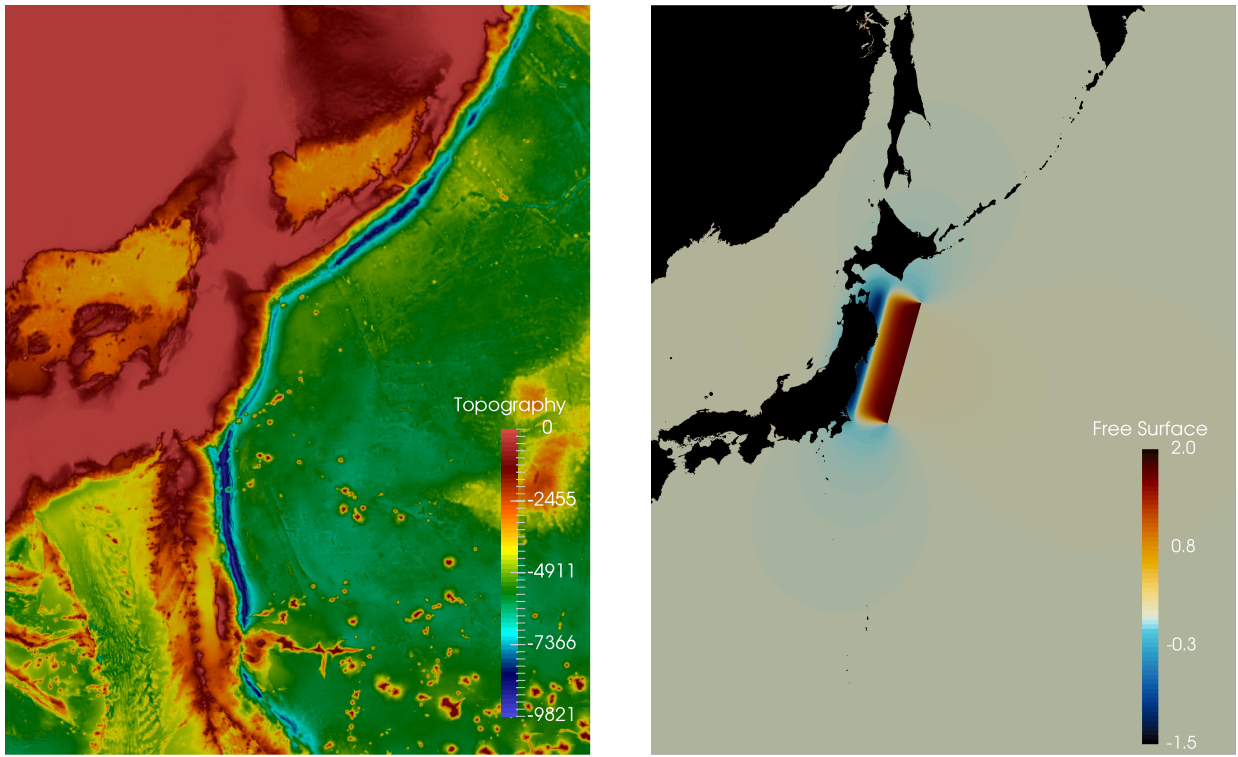

Figure 18: Depiction of the 2011 tsunami simulation. Left panel: bathymetry (submerged topography). The Kuril trench is the deepest part of the ocean, depicted in deep blue. The continents are represented in red. Right panel: initial free surface. The continents are depicted in black, the average water surface in gray, and the initial tsunami wave lies over the Kuril trench, next to the Japanese coast.

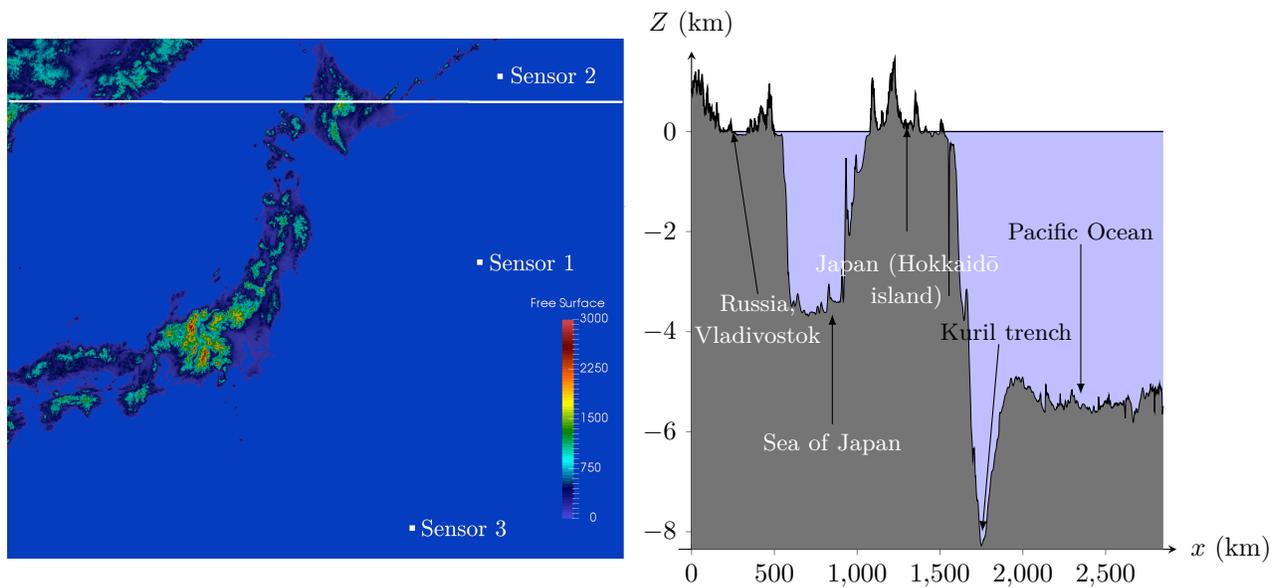

Figure 19: Left panel: position of the three sensors. Right panel: depiction of the topography over the white line drawn in the left panel.

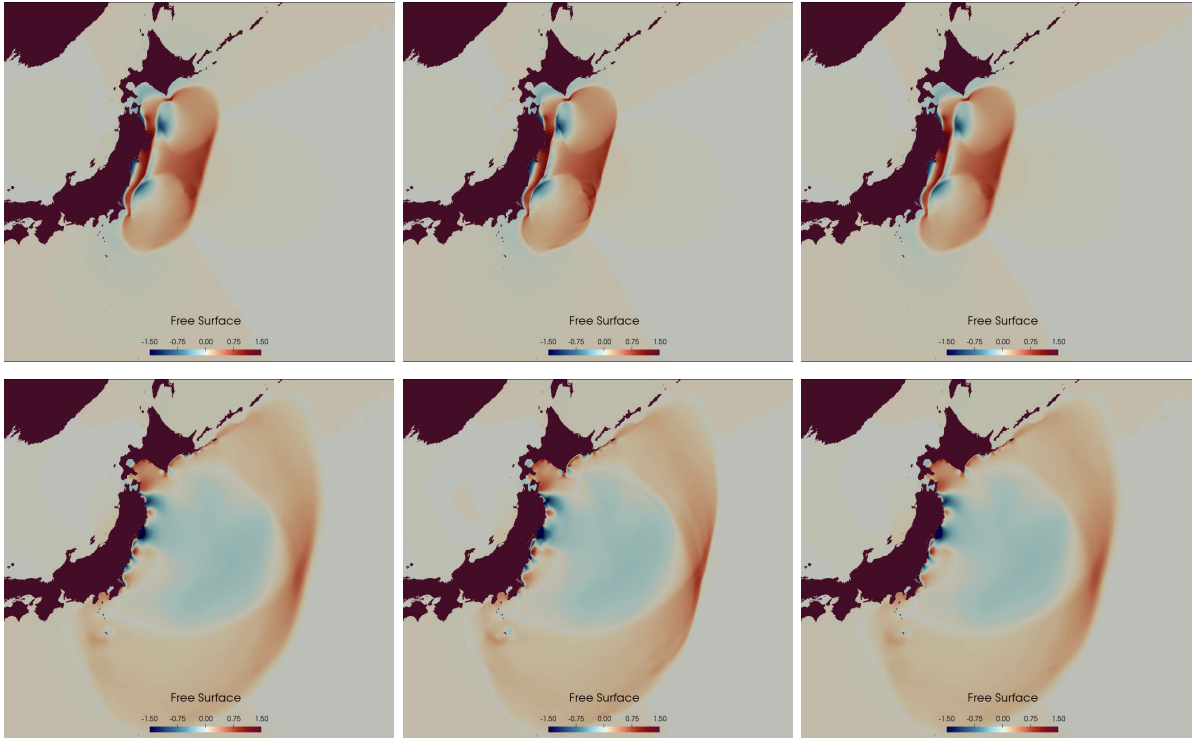

Figure 20: Simulation of the 2011 Japan tsunami with the $\mathbb{P}_0$ scheme (left panels), the $\mathbb{P}_1^{\text{WB}}$ scheme (middle panels) and the $\mathbb{P}_3^{\text{WB}}$ scheme (right panels). Snapshots taken at times $t = 720$s (top panels) and $t = 3600$s (bottom panels). The average sea surface height is represented in gray, and the continents are displayed in dark red.

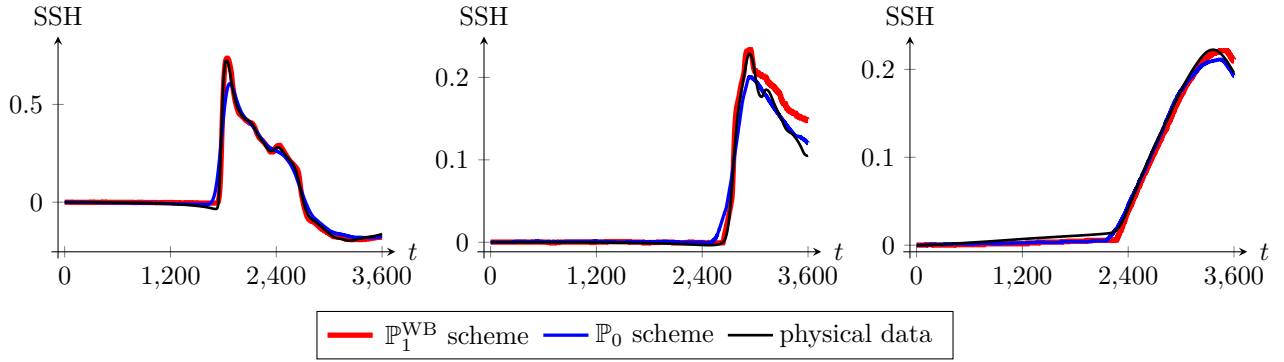

Figure 21: Sea surface height at each sensor (from left to right, sensor #1 to sensor #3, whose positions are displayed in Figure 19). The physical data is represented in black, the $\mathbb{P}_0$ approximation in blue, and the $\mathbb{P}_1^{\text{WB}}$ approximation in red. The total water depths below sensors #1, #2 and #3 are respectively 5700m, 6600m and 4400m.

*Acknowledgments*

V. Michel-Dansac extends his thanks to the Service Hydrographique et Océanographique de la Marine (SHOM) for financial support. C. Berthon and F. Foucher would like to thank the project MUFFIN ANR-19-CE46-0004 for financial support. S. Clain acknowledges the financial support by FEDER – Fundo Europeu de Desenvolvimento Regional, Portugal, through COMPETE 2020 – Programa Operational Fatores de Competitividade, and the National Funds through FCT – Fundação para a Ciência e a Tecnologia, Portugal, project No. UID/FIS/04650/2013 and project No. POCI-01-0145-FEDER-028118.